\documentclass[12pt]{article}
\usepackage{latexsym,amsmath}
\usepackage{amsfonts}
\usepackage{amsmath}
\usepackage{amssymb}
\usepackage{latexsym}
\usepackage{indentfirst}
\usepackage{fancyhdr}
\usepackage{stmaryrd}
\usepackage{mathrsfs}
\usepackage{cite}
\date{}

\numberwithin{equation}{section}

\newtheorem{theorem}{Theorem}
\newtheorem{lemma}[theorem]{Lemma}

\begin{document}

\newcommand{\for}{\mbox{for}}
\newcommand{\ptl}{\partial}

\title{Twisted Hamiltonian Lie Algebras and Their Multiplicity-Free
Representations }
\author{Ling Chen
\thanks{Corresponding author. E-mail: chenling@amss.ac.cn (L. Chen)}\\%EndAName
{\small{ Institute of Applied Mathematics, Academy of Mathematics \& System Sciences,}}\\
{\small{ Chinese Academy of Sciences, Beijing 100190,
China}}}\maketitle

\begin{abstract}

We construct a class of new Lie algebras by generalizing the
one-variable Lie algebras generated by the quadratic conformal
algebras (or corresponding Hamiltonian operators) associated to
Poisson algebras and a quasi-derivation found by Xu. These algebras
can be viewed as certain twists of Xu's generalized Hamiltonian Lie
algebras. The simplicity of these algebras is completely determined.
Moreover, we construct a family of multiplicity-free representations
of these Lie algebras and prove their irreducibility.

\end{abstract}
\section{Introduction}

Simple Lie algebras of Cartan type play important roles in the
structure theory of Lie algebras. They also have many applications
in the areas of geometry and mathematical physics. After Kac
\cite{Ka} gave an abstract definition of Lie algebras of Cartan type
by derivations, various graded generalizations were introduced and
studied by Kawamoto \cite{K}, Osborn \cite{O}, Dokovic and Zhao
\cite{DZ}, and Zhao \cite{Zk1}. Motivated from his works on
quadratic conformal algebras \cite{X3} and Hamiltonian operators
\cite{X4}, Xu \cite{X2} found certain nongraded generalizations of
the Lie algebras of Cartan type and determined their simplicity.
Indeed, Xu's algebras in general do not contain any toral Cartan
subalgebras, so traditional approaches to simple Lie algebras do not
work in his case. In particular, Xu's generalizations of Hamiltonian
Lie algebras were constructed from a skew-symmetric
$\mathbb{Z}$-bilinear form and an even number of locally-finite
commuting derivations. Hamiltonian Lie algebras play key roles in
both classical and quantum physics. They are also fundamental
algebraic structure in symplectic geometry.

Passman \cite{P} gave a certain necessary and sufficient condition
on derivations for a Lie algebra of generalized Witt type to be
simple. The isomorphic structural spaces of Xu's nongrared Lie
algebras of Cartan type were determined in \cite{SX1, SX2, SX3,
SXZ}. Moreover, Xu \cite{X1} also introduced certain nongraded
generalizations of the Lie algebras of Block type and determined
their simplicity.

On the representation side,  Shen \cite{Sg1, Sg2, Sg3} introduced
mixed product of graded modules over graded Lie algebras of Cartan
type (also known as Larsson functor (cf. \cite{La})) and obtained
certain irreducible modules. Rao \cite{R1, R2} constructed some
irreducible weight modules over the derivation Lie algebra of the
algebra of Laurent polynomials based on Shen's mixed product. Lin
and Tan \cite{LT} did the similar thing over the derivation Lie
algebra of the algebra of quantum torus.  Zhao \cite{Zy} determined
the module structure of Shen's mixed product over Xu's nongraded Lie
algebras of Witt type.

Multiplicity-free representations are important objects in Lie
theory due to their  neat formulae and natural applications. Su and
Zhao \cite{SZk} classified the multiplicity-free representations
over generalized Virasoro algebras. Moreover, Zhao \cite{Zk2} gave a
classification of the multiplicity-free representations over graded
generalized Witt algebras.  Furthermore, Su and Zhou \cite{SZ}
classified the analogues of multiplicity-free representations over
Xu's nongraded generalized Witt algebras.

In \cite{X3}, Xu classified certain quadratic conformal algebras and
constructed such algebras from Poisson algebras and a
quasi-derivation. Conformal algebras are the local structures of
certain one-variable Lie algebras. They are also equivalent to
certain linear Hamiltonian operators (cf. \cite{X4}). We construct a
class of new Lie algebras by generalizing the one-variable Lie
algebras generated by the quadratic conformal algebras associated to
Poisson algebras and a quasi-derivation. These algebras can be
viewed as certain twists of Xu's generalized Hamiltonian Lie
algebras. The simplicity of these algebras is completely determined.
Moreover, we construct a family of multiplicity-free representations
of these Lie algebras and prove their irreducibility. Below we give
more technical details.

Throughout this paper, let $\mathbb{F}$ be a field with
characteristic 0. All the vector spaces are assumed over
$\mathbb{F}$. Let $({\cal A},[\cdot,\cdot],\cdot)$ be a Poisson
algebra, that is, $({\cal A},[\cdot,\cdot])$ forms a Lie algebra,
$({\cal A},\cdot)$ forms a commutative associative algebra and the
following compatibility condition holds:
\begin{equation}
[u,v\cdot
w]=[u,v]\cdot w+v\cdot [u,w]\qquad\for\;u,v,w\in{\cal
A}.
\end{equation}
Let $\ptl_1$ and $\ptl_2$ be two commuting
derivations of $({\cal A},\cdot)$. Suppose that $\ptl_2$ is also a
derivation of the Lie algebra $({\cal A},[\cdot,\cdot])$ but
$\ptl_1$ is a quasi-derivation of $({\cal A},[\cdot,\cdot])$ in the
following sense:
\begin{equation}
\ptl_1([u,v])=[\ptl_1(u),v]+[u,\ptl_1(v)]+c[u,v]\qquad\for\;u,v\in{\cal
A},
\end{equation}
where $c\in\mathbb{F}$ is a fixed constant. Define
a new algebraic operation $[\cdot,\cdot]_1$ on ${\cal A}$ by
\begin{equation}\label{equ:5.1}
[u,v]_1=[u,v]+(\ptl_1+c)(u)\cdot\ptl_2(v)-\ptl_2(u)\cdot(\ptl_1+c)(v)\qquad\for\;u,v\in{\cal
A}.
\end{equation}
It can be verified that $({\mathcal A},[\cdot,\cdot]_1)$ also forms
a Lie algebra, which is motivated from the one-variable Lie algebras
generated by the quadratic conformal algebras associated to Poisson
algebras and a quasi-derivation in \cite{X3}. If $[\cdot,\cdot]$ is
the trivial Lie bracket, that is, $[u,v]=0$ for $u,v\in{\cal A}$,
then $({\cal A},[\cdot,\cdot]_1)$ with $c=1$ is exactly Xu's
generalization of Block algebra in \cite{X1}. In the case $c=0$,
(\ref{equ:5.1}) can be viewed as the iterated construction of
Hamiltonian Lie algebras. When ${\cal
A}=\mathbb{F}[x_1,...,x_{2n+1}]$, $\ptl_2=\ptl_{x_{2n+1}}$, $\ptl_1$
is the grading operator on ${\cal A}$ with respect to $\{x_1
,...,x_{2n}\}$ and
\begin{equation}[f,g]=\sum_{i=1}^n(\ptl_{x_i}(f)\ptl_{x_{n+i}}(g)-\ptl_{x_{n+i}}(f)\ptl_{x_i}(g))
\qquad\for\;f,g\in {\cal A},\end{equation} the Lie algebra $({\cal
A},[\cdot,\cdot]_1)$ with $c=-2$ is exactly the canonical Lie
algebra of Contact type. In this paper, we will deal with another
case, which can be viewed as a generalization of Xu's
four-derivation analogues of Block algebras in \cite{X1}.

Denote by $\mathbb{Z}$ the ring of integers and by $\mathbb{N}$ the
additive semigroup of nonnegative integers.  For two integers $i\leq
j$, we denote $\overline{i,j}=\{i,i+1,\cdots,j\}$. Let $n$ be a
positive integer. Picking
\begin{equation}
\mathcal {J}_{p}\in\{\{0\},\mathbb{N}\}\qquad \textrm{for}\
p\in\overline{1,2n},
\end{equation}
we have the following product of
additive semigroups
\begin{equation}
\boldsymbol{\mathcal{J}}=\mathcal{J}_{1}\times\mathcal{J}_{2}\times\cdots\times\mathcal{J}_{2n}.
\end{equation}

Let $\Gamma$ be a nonzero torsion-free abelian group, which is also
viewed as a $\mathbb{Z}$-module throughout this paper. Let
$\phi(\cdot,\,\cdot): \Gamma \times \Gamma \rightarrow \mathbb{F}$
be a skew-symmetric $\mathbb{Z}$-bilinear form and let
$\{\varphi_{p}: \Gamma \rightarrow \mathbb{F} \mid
p\in\overline{1,2n}\}$  be $2n$ additive group homomorphisms such
that
\begin{equation}\label{equ:1.1}
\varphi_{p}(\Gamma)+\mathcal{J}_{p}\not=\{0\} \qquad\textrm{for}\
p\in\overline{1,2n},
\end{equation}
\begin{equation}\label{equ:1.3}
\textrm{Rad}_{\phi} \bigcap \left( \bigcap _{ p \in \overline{1,2n}}
\textrm{ker}_{\varphi_p} \right) = \{ 0 \},
\end{equation}
and
\begin{equation}\label{equ:1.2}
\textrm{Rad}_{\phi} \bigcap \left( \bigcap _{q \not= p \in
\overline{1,2n}} \textrm{ker}_{\varphi_p} \right) \backslash
\textrm{ker}_{\varphi_q}\not = \emptyset , \ \ \textrm{if }
\varphi_{q} \not \equiv 0, \ \ q\in \overline{1,2n}.
\end{equation}

 We choose
fixed elements
\begin{equation}\label{equ:1.4}
\sigma_{2s+1} = \sigma_{2s+2} \in \textrm{Rad}_{\phi} \bigcap \left(
\bigcap _{2s+1, 2s+2 \not= p \in \overline{1,2n}}
\textrm{ker}_{\varphi_p} \right)
\end{equation}
such that
$\sigma_{2s+1}\not=0$ if $\varphi_{2s+1} \not\equiv 0$ or
$\varphi_{2s+2} \not\equiv 0$ for $s\in \overline{1,n-1}$.

Let $\mathscr{A}$ be a vector space with a basis $
\{x^{\alpha,\mathbf{i}}\mid(\alpha,\mathbf{i})\in\Gamma\times\boldsymbol{\mathcal{J}}\}.$
Define an algebraic operation $\cdot$ on $\mathscr{A}$ by
\begin{equation}
x^{\alpha,\mathbf{i}}\cdot
x^{\beta,\mathbf{j}}=x^{\alpha+\beta,\mathbf{i}+\mathbf{j}}\qquad\textrm{for}\
(\alpha,\mathbf{i}),(\beta,\mathbf{j})\in\Gamma\times\boldsymbol{\mathcal{J}}.\end{equation}
Then $\mathscr{A}$ forms a commutative associative algebra with an
identity element $x^{0,\mathbf{0}}$, which will simply be denoted by
$1$ in the rest of this paper. Throughout this paper, the notion
``$\cdot$'' will be invisible in a product when the context is
clear.

Use $\mathbf{i} =(i_1, i_2, \cdots, i_{2n})$ to denote an element in
$\boldsymbol{\mathcal{J}}$. Moreover, we denote
\begin{equation}
a_{[p]}=(0,\cdots,\stackrel{p}{a},0,\cdots,0),\qquad
a\in\mathcal{J}_{p}\,.\end{equation} For $p\in\overline{1,2n}$, we
define $\partial_p\in\textrm{End}\mathscr{A}$ by
\begin{equation} \label{equ:1.6}
\partial_p(x^{\alpha,\mathbf{i}})=\varphi_p(\alpha)x^{\alpha,\mathbf{i}}+i_p x^{\alpha,\mathbf{i}-1_{[p]}}
\end{equation}
for $(\alpha,\mathbf{i})\in\Gamma\times\boldsymbol{\mathcal{J}}$,
where we adopt the convention that if a notion is not defined but
technically appears in an expression, we always treat it as zero;
for instance, $x^{\alpha,-1_{[1]}}=0$ for any $\alpha\in\Gamma$. It
can be verified that $\{\partial_p\mid p\in \overline{1,2n} \}$ are
commuting derivations of $\mathscr{A}$. For any $\alpha \in \Gamma$,
we set $\mathscr{A}_{\alpha} = \textrm{Span} \{
x^{\alpha,\mathbf{j}}\mid \mathbf{j} \in
\boldsymbol{\mathcal{J}}\}.$ Assume that
$\textrm{Rad}_{\phi}\bigcap(\bigcap_{p=2}^{2n}\textrm{ker}_{\varphi_p})$
has an element $\varepsilon$ such that $\varphi_{1}(\varepsilon)=1$.
Then we have a Lie bracket $[\cdot, \cdot]$ of the type
(\ref{equ:5.1}) on $\mathscr{A}$ determined by
\begin{eqnarray}\label{equ:1.8}
[u,w] & = & x^{\varepsilon, \boldsymbol{0}} \left( \sum^{n-1}_{s=1}
x^{\sigma_{2s+1}, \boldsymbol{0}}
(\partial_{2s+1}(u)\partial_{2s+2}(w)-\partial_{2s+1}(w)\partial_{2s+2}(u))
+ \phi(\alpha, \beta)uw \right) \nonumber\\
& &
 +(\partial_{1}+1)(u)\partial_{2}(w)-(\partial_{1}+1)(w)\partial_{2}(u)
\end{eqnarray}
for $ u \in \mathscr{A}_{\alpha}$ and $w \in \mathscr{A}_{\beta}$.

Notice that $x^{-\varepsilon, \boldsymbol{0}}$ is a central element
of $\mathscr{A}$. Form a quotient Lie algebra
$\mathscr{H}=\mathscr{A}/\mathbb{F}x^{-\varepsilon,\boldsymbol{0}}$,
whose induced Lie bracket is still denoted by $[\cdot,\cdot]$ when
the context is clear. We call the Lie algebra $(\mathscr{H},
[\cdot,\cdot])$ a {\it twisted Hamiltonian Lie algebra}. Denote
$\sigma=\sum_{s=1}^{n-1} \sigma_{2s+1}-2\varepsilon$.

\begin{theorem}\label{theorem:1}
The Lie algebra ($\mathscr{H},[\cdot,\cdot]$) is simple if
$\boldsymbol{\mathcal{J}}\not=\{\boldsymbol{0}\}$. If
$\boldsymbol{\mathcal{J}}=\{\boldsymbol{0}\}$, then
$\mathscr{H}^{(1)}=[\mathscr{H},\mathscr{H}]$ is a simple Lie
algebra and $\mathscr{H}=\mathscr{H}^{(1)}\oplus (\mathbb{F}
x^{\sigma, \boldsymbol{0}} + \mathbb{F}x^{-\varepsilon,
\boldsymbol{0}})$.
\end{theorem}

Take $\vec{\xi}=(\xi_{1}, \cdots, \xi_{2n}) \in \mathbb{F}^{2n}$ and
an additive group homomorphism $f:\Gamma\rightarrow\mathbb{F}$ such
that $\textrm{Rad}_{\phi} \subset \textrm{ker}f$. For $u\in
\mathscr{A}$, we denote $\bar u=u+\mathbb{F}x^{-\varepsilon,
\boldsymbol{0}}$. In particular, we denote
$\bar{x}^{\alpha,\mathbf{i}}=x^{\alpha,\mathbf{i}}+\mathbb{F}x^{-\varepsilon,
\boldsymbol{0}}$ for
$(\alpha,\mathbf{i})\in\Gamma\times\boldsymbol{\mathcal{J}}$. It can
be verified that $\mathscr{A}$ forms a {\it multiplicity-free}
$\mathscr{H}$-module with the action determined by
\begin{eqnarray}\label{equ:1.9}
\bar u.v & =&  x^{\varepsilon, \boldsymbol{0}} \big(
\sum^{n-1}_{s=1} x^{\sigma_{2s+1}, \boldsymbol{0}}
(\partial_{2s+1}(u)(\partial_{2s+2}+\xi_{2s+2})(v)-\partial_{2s+2}(u)(\partial_{2s+1}+\xi_{2s+1})(v))
\nonumber\\& +&(\phi(\alpha, \beta)+f(\alpha))uv \big)
 +(\partial_{1}+1)(u)(\partial_{2}+\xi_{2})(v)-\partial_{2}(u)(\partial_{1}+\xi_1)(v)
\end{eqnarray}
for $ u \in \mathscr{A}_{\alpha}$ and $v \in \mathscr{A}_{\beta}$.
 To make it distinguishable,
we denote by $\mathscr{A}_{\vec{\xi},f}$ the module $\mathscr{A}$
with the above action. In particular, denoting $\vec{\xi_{0}}=(0,
\cdots, 0, 1)\in \mathbb{F}^{2n}$, module $\mathscr{A}_{\vec
{\xi_{0}},0}$ is induced from the adjoint representation of the Lie
algebra $\mathscr{A}$. Moreover, it is straightforward to verify
that $\mathscr{A}_{\vec{\xi}',f'} \cong \mathscr{A}_{\vec{\xi},f}$
if there exists $\mu' \in \Gamma$ such that
$\xi'_{p}-\xi_{p}=\varphi_{p}(\mu')\textrm{ for }
p\in\overline{1,2n},$ and $ f'(\beta)-f(\beta) = \phi(\beta,\mu')
\textrm{ for any } \beta \in \Gamma$.

\begin{theorem}\label{theorem:2}
If there exists $\mu \in \Gamma$ such that
\begin{equation}\label{equ:1.10}
\xi_{p}=\varphi_{p}(\mu) \quad \textrm{ for } p\in\overline{1,2n},
\end{equation}
and
\begin{equation}\label{equ:1.11}
f(\alpha) = \phi(\alpha,\mu) \quad \textrm{ for any } \alpha \in
\Gamma,
\end{equation}
then $\mathscr{A}_{\vec{\xi},f}\cong \mathscr{A}_{\vec{\xi_{0}},0}$.
Otherwise, $\mathscr{A}_{\vec{\xi},f}$ is irreducible.
\end{theorem}

We shall present the proofs of Theorem 1 and Theorem 2 in sections 2
and 3, respectively.

\section{Proof of Theorem 1}

In this section, we will determine the simplicity of the twisted
Hamitonian Lie algebra $\mathscr{H}$. In order to prove simplicity,
we need the following Lemma.

\begin{lemma}\label{lemma:1} Let $\mathbf{T}$ be a linear transformation on a vector space $U$,
and let $U_1$ be a subspace of $U$ such that $\mathbf{T}
(U_1)\subset U_1$. Suppose that $u_1,\ u_2,\ \cdots,\ u_n$ are
eigenvectors of $\mathbf{T}$ corresponding to different eigenvalues.
If $\sum_{j=1}^n u_j \in U_1$, then $u_1,\ u_2,\ \cdots,\ u_n\ \in
U_1$.
\end{lemma}

Let $\mathcal {I}$ be any ideal of $\mathscr{A}$ that strictly
contains $\mathbb{F}x^{-\varepsilon, \boldsymbol{0}}$. When
$\boldsymbol{\mathcal{J}}\not=\{\boldsymbol{0}\}$, proving the
simplicity of ($\mathscr{H},[\cdot,\cdot]$) is equivalent to proving
that $\mathcal {I}=\mathscr{A}$.

Our proof of Theorem \ref{theorem:1} is organized as follows.
Firstly, we will determine some preliminary facts about the ideal
$\mathcal {I}$. Then, we will process the proof in four steps.

First by (\ref{equ:1.8}), we have
\begin{eqnarray}\label{equ:2.1}
& & [x^{\alpha,\mathbf{i}}, x^{\beta,\mathbf{j}}]\nonumber\\
 & = &((\varphi_1(\alpha)+1) \varphi_2(\beta)-\varphi_2(\alpha)
(\varphi_1(\beta)+ 1))x^{\alpha+\beta,\mathbf{i}+\mathbf{j}}+ (i_{1}
\varphi_2(\beta)-j_{1} \varphi_2(\alpha))\cdot \nonumber\\ & &
x^{\alpha+\beta,\mathbf{i}+\mathbf{j}-1_{[1]}} + (j_{2}
(\varphi_1(\alpha)+1) - i_{2} (\varphi_1(\beta)+1))
x^{\alpha+\beta,\mathbf{i}+\mathbf{j} -1_{[2]}}
 + (i_{1}j_{2}-i_{2}j_{1})\cdot\nonumber\\& & x^{\alpha+\beta,\mathbf{i}+\mathbf{j}-1_{[1]}-1_{[2]}}
  + x^{\varepsilon,\boldsymbol{0}}\Big(\sum^{n-1}_{s=1} x^{\sigma_{2s+1},\boldsymbol{0}}
  \big((\varphi_{2s+1}(\alpha) \varphi_{2s+2}(\beta)-\varphi_{2s+2}(\alpha)
\varphi_{2s+1}(\beta))\nonumber\\
 & &\cdot
x^{\alpha+\beta,\mathbf{i}+\mathbf{j}} +(i_{2s+1}
\varphi_{2s+2}(\beta)-j_{2s+1}
\varphi_{2s+2}(\alpha))x^{\alpha+\beta,\mathbf{i}+\mathbf{j}-1_{[2s+1]}}
+(j_{2s+2} \varphi_{2s+1}(\alpha) \nonumber\\
 & &- i_{2s+2}
\varphi_{2s+1}(\beta))x^{\alpha+\beta,\mathbf{i}+\mathbf{j}-1_{[2s+2]}}
 +(i_{2s+1}j_{2s+2}-i_{2s+2}j_{2s+1})\cdot\nonumber\\
  & &  x^{\alpha+\beta,\mathbf{i}+\mathbf{j}-1_{[2s+1]}-1_{[2s+2]}}\big)+\phi(\alpha,\beta)x^{\alpha+\beta,\mathbf{i}+\mathbf{j}}
\Big)\end{eqnarray} for
$(\alpha,\mathbf{i}),(\beta,\mathbf{j})\in\Gamma\times\boldsymbol{\mathcal{J}}$.
In particular,
\begin{equation}\label{equ:2.2}
[1, x^{\beta, \mathbf{j}}]=\varphi_{2}(\beta)x^{\beta,
\mathbf{j}}+j_{2}x^{\beta, \mathbf{j}-1_{[2]}} \qquad\textrm{for}\
(\beta,\mathbf{j})\in\Gamma\times\boldsymbol{\mathcal{J}}.
\end{equation}
Set
\begin{equation}\label{equ:2.3}
\Gamma'=\{\varphi_{2}(\alpha)\mid \alpha\in\Gamma\}.
\end{equation}
For $\lambda\in\Gamma'$, we define
\begin{equation}\label{equ:2.4}
\mathscr{A}_{(\lambda)}=\{u\in\mathscr{A} \mid
(\textrm{ad}_{1}-\lambda)^{m}(u)=0\ \textrm{for some}\
m\in\mathbb{N}\}.
\end{equation}
Then by ($\ref{equ:2.2}$), we have
\begin{equation}\label{equ:2.5}
\mathscr{A}=\bigoplus_{\lambda\in\Gamma'}\mathscr{A}_{(\lambda)}.
\end{equation}
In particular,
\begin{equation}\label{equ:2.6}
\mathcal{I}=\bigoplus_{\lambda\in\Gamma'}\mathcal{I}_{\lambda},\qquad
\mathcal{I}_{\lambda}=\mathcal{I}\cap\mathscr{A}_{(\lambda)}.
\end{equation}
For any $(\beta,\mathbf{j})\in\Gamma\times\boldsymbol{\mathcal{J}}$,
we find
\begin{eqnarray}
[x^{-\beta, \mathbf{0}}, x^{\beta, \mathbf{j}}] & =&
\sum^{n-1}_{s=1}x^{\varepsilon+\sigma_{2s+1},\boldsymbol{0}}(j_{2s+1}
   \varphi_{2s+2}(\beta)x^{0,\mathbf{j}-1_{[2s+1]}}
 -j_{2s+2}\varphi_{2s+1}(\beta)x^{0,\mathbf{j}-1_{[2s+2]}})\nonumber\\
   &&+2\varphi_2(\beta)x^{0,\mathbf{j}}
 +j_{1} \varphi_2(\beta)x^{0,\mathbf{j}-1_{[1]}}
  - j_{2} (\varphi_1(\beta)-1)x^{0,\mathbf{j}-1_{[2]}} \label{equ:2.7}
\end{eqnarray}
by (\ref{equ:2.1}). If $\mathcal{I}_{\lambda}\not=\{0\}$ for some
$0\not=\lambda\in\Gamma'$, we pick any
\begin{equation}0\not=u=\sum_{(\alpha,\mathbf{i})\in\Gamma\times\boldsymbol{\mathcal{J}}}d_{\alpha,\mathbf{i}}
x^{\alpha,\mathbf{i}}\in\mathcal{I}_{\lambda}, \qquad
d_{\alpha,\mathbf{i}}\in\mathbb{F}.\end{equation} Assume
$d_{\beta,\mathbf{j}}\not=0$ for some
$(\beta,\mathbf{j})\in\Gamma\times\boldsymbol{\mathcal{J}}$. Denote
\begin{equation}|\mathbf{i}| =\sum_{p=1}^{2n} i_{p}\qquad\for\;\mathbf{i}\in \boldsymbol{\mathcal{J}}.
\end{equation}
 Fixing
$\beta$ and considering the nonzero terms
$d_{\beta,\mathbf{l}}x^{\beta,\mathbf{l}}$ in $u$ with the largest
value $|\mathbf{l}|$, we have
\begin{equation}\label{equ:2.8}
[x^{-\beta, \mathbf{0}}, u]\in
\mathcal{I}_{0}\backslash\mathbb{F}x^{-\varepsilon,\mathbf{0}}
\end{equation}
by (\ref{equ:2.7}) and the fact  $\varepsilon\not=\mathbf{0}$. Thus
(\ref{equ:2.6}) implies that $\mathcal{I}_{0}$ strictly contains
$\mathbb{F}x^{-\varepsilon,\mathbf{0}}$.

Observe
\begin{equation}\label{equ:2.9}
[1, x^{\beta, \mathbf{j}}]=j_{2}x^{\beta, \mathbf{j}-1_{[2]}}
\qquad\textrm{for}\
(\beta,\mathbf{j})\in\textrm{ker}_{\varphi_{2}}\times\boldsymbol{\mathcal{J}}
\end{equation}
and
\begin{equation}\label{equ:2.10}
[x^{\varepsilon,\mathbf{0}},x^{-\varepsilon,1_{[2]}}]=2x^{0,\mathbf{0}}
\end{equation}
by (\ref{equ:2.1}). Set
\begin{equation}
\mathcal{I}_{[0]}=\mathcal{I}\cap (\sum_{(\alpha, \mathbf{i})\in
\textrm{ker}_{\varphi_{2}} \times \boldsymbol{\mathcal{J}};\:
i_{2}=0} \mathbb{F} x^{\alpha, \mathbf{i}}).
\end{equation}
By (\ref{equ:2.9}) and (\ref{equ:2.10}) and repeatedly acting of
$\textrm{ad}_{1}$ on $\mathcal{I}_{0}$, we can prove that
\begin{equation}
\mathcal{I}_{[0]}\ \textrm{strictly contains}\ \mathbb{F}
x^{-\varepsilon, \mathbf{0}}.
\end{equation}

For $ k \in \mathbb{N} $, let
\begin{equation}
\mathscr{A}_{\langle k \rangle}=\textrm{span}
\{x^{\alpha,\mathbf{i}} \mid
(\alpha,\mathbf{i})\in\Gamma\times\boldsymbol{\mathcal{J}},|\mathbf{i}|\leq
k \}.
\end{equation}
Then we have
\begin{equation}
\mathscr{A}=\bigcup_{k=0}^{\infty} \mathscr{A}_{\langle k \rangle}.
\end{equation}
Set
\begin{equation}\label{equ:2.11}
\hat{k}=\textrm{min} \{ k \in \mathbb{N} \mid (\mathcal{I}_{[0]}\cap
\mathscr{A}_{\langle k \rangle}) \backslash
\mathbb{F}x^{-\varepsilon, \mathbf{0}} \not = \emptyset \}.
\end{equation}
For any $u \in (\mathcal{I}_{[0]}\cap \mathscr{A}_{\langle \hat{k}
\rangle}) \backslash \mathbb{F}x^{-\varepsilon, \mathbf{0}}$, we
write it as
\begin{equation}\label{equ:2.12}
u= \sum_{\substack{(-\varepsilon, \mathbf{0})\not =(\alpha,
\mathbf{i})\in \textrm{ker}_{\varphi_{2}} \times
\boldsymbol{\mathcal{J}}\\ i_{2}=0, |\mathbf{i}|= \hat{k}}}
a_{\alpha, \mathbf{i}} x^{\alpha, \mathbf{i}}+u',
\end{equation}
where $a_{\alpha, \mathbf{i}} \in \mathbb{F},\ u' \in
\mathscr{A}_{\langle \hat{k}-1 \rangle}+\mathbb{F}x^{-\varepsilon,
\mathbf{0}}$. Moreover, we define
\begin{equation}\label{equ:2.13}
\iota (u)=|\{\alpha\in\textrm{ker}_{\varphi_{2}}\mid a_{\alpha,
\mathbf{i}} \not= 0 \ \textrm{for some}\ \mathbf{i}\in
\boldsymbol{\mathcal{J}}, |\mathbf{i}|= \hat{k} \}|.
\end{equation}
By (\ref{equ:2.11}), $\iota (u)>0$. Furthermore, we set
\begin{equation}\label{equ:2.14}
\iota = \textrm{min} \{ \iota (w) \mid w \in (\mathcal{I}_{[0]}\cap
\mathscr{A}_{\langle \hat{k} \rangle}) \backslash
\mathbb{F}x^{-\varepsilon, \mathbf{0}} \}.
\end{equation}
Choose $u \in (\mathcal{I}_{[0]}\cap \mathscr{A}_{\langle \hat{k}
\rangle}) \backslash \mathbb{F}x^{-\varepsilon, \mathbf{0}} $ such
that $\iota (u) = \iota$. Write $u$ as in (\ref{equ:2.12}). We will
process our proof in four steps.\vspace{0.3cm}

\noindent{\bf {\emph{Step 1. $\alpha=\beta$ whenever $\ a_{\alpha,
\mathbf{i}}a_{\beta, \mathbf{j}} \not= 0.$}}} \vspace{0.3cm}

Firstly, we want to prove that
\begin{equation}\label{equ:2.15}
\varphi_{1}(\alpha)=\varphi_{1}(\beta)\ \textrm{whenever}\
a_{\alpha, \mathbf{i}}a_{\beta, \mathbf{j}} \not= 0.
\end{equation}

If $\mathcal{J}_{2} = \mathbb{N}$, we get
\begin{eqnarray}\label{equ:2.16}
&&\mathcal{I}_{[0]}\cap \mathscr{A}_{\langle \hat{k} \rangle} \owns
[x^{0, 1_{[2]}}, u]\nonumber\\
 & \equiv & \sum_{\substack{(-\varepsilon, \mathbf{0})\not =(\alpha,
\mathbf{i})\in
\textrm{ker}_{\varphi_{2}} \times \boldsymbol{\mathcal{J}}\\
i_{2}=0, |\mathbf{i}|= \hat{k}}} -a_{\alpha, \mathbf{i}}
(\varphi_{1}(\alpha)+1) x^{\alpha, \mathbf{i}}\;\; (\textrm{mod}\
\mathscr{A}_{\langle \hat{k}-1 \rangle}).
\end{eqnarray}
By the minimality of $\iota (u)$ and Lemma \ref{lemma:1}, we have
\begin{equation}\label{equ:2.17}
\varphi_{1}(\alpha)+1 = \varphi_{1}(\beta)+1 \quad
\textrm{whenever}\ a_{\alpha, \mathbf{i}}a_{\beta, \mathbf{j}} \not=
0,
\end{equation}
that is,
\begin{equation}\label{equ:2.18}
\varphi_{1}(\alpha)=\varphi_{1}(\beta)\  \quad \textrm{whenever}\
a_{\alpha, \mathbf{i}}a_{\beta, \mathbf{j}} \not= 0.
\end{equation}
Then (\ref{equ:2.15}) holds.

If $\mathcal{J}_{2} = \{0\}$, then $\varphi_{2} \not \equiv 0$ by
(1.7). Picking $\tau \in \textrm{Rad}_{\phi} \cap (\cap_{2\not=p \in
\overline{1,2n}} \textrm{ker}_{\varphi_{p}}) \backslash
\textrm{ker}_{\varphi_{2}}$, we have
\begin{eqnarray}\label{equ:2.19}
\mathcal{I}_{[0]}\cap \mathscr{A}_{\langle \hat{k} \rangle} &\owns&
[x^{-\tau, \mathbf{0}},[x^{\tau, \mathbf{0}}, u]] \equiv
\sum_{\substack{(-\varepsilon, \mathbf{0})\not =(\alpha,
\mathbf{i})\in
\textrm{ker}_{\varphi_{2}} \times \boldsymbol{\mathcal{J}}\\
i_{2}=0, |\mathbf{i}|= \hat{k}}} -a_{\alpha, \mathbf{i}}
\varphi^{2}_{2}(\tau) (\varphi_{1}(\alpha)+1)\cdot\nonumber\\
 && (\varphi_{1}(\alpha)+2)
x^{\alpha, \mathbf{i}}\;\;(\textrm{mod}\ \mathscr{A}_{\langle
\hat{k}-1 \rangle}).
\end{eqnarray}
By the minimality of $\iota (u)$ and Lemma \ref{lemma:1}, we get
\begin{equation}\label{equ:2.20}
(\varphi_{1}(\alpha)+1)(\varphi_{1}(\alpha)+2) =
(\varphi_{1}(\beta)+1)(\varphi_{1}(\beta)+2) \
\end{equation}
whenever $\ a_{\alpha, \mathbf{i}}a_{\beta, \mathbf{j}} \not= 0$,
which implies
\begin{equation}\label{equ:2.21}
\varphi_{1}(\alpha)=\varphi_{1}(\beta)\ \textrm{ or }
\varphi_{1}(\alpha+\beta)+3=0 \ \quad \textrm{whenever}\ a_{\alpha,
\mathbf{i}}a_{\beta, \mathbf{j}} \not= 0.
\end{equation}
Assume that there exist $\alpha, \beta \in \Gamma$ such that
$\varphi_{1}(\alpha) \not= \varphi_{1}(\beta),\
\varphi_{1}(\alpha+\beta)+3=0$ and $a_{\alpha, \mathbf{i}}a_{\beta,
\mathbf{j}} \not= 0$. Without loss of generality, we assume that
$\varphi_{1}(\alpha)+1\not= 0$. Note that $\varphi_{1} \not \equiv
0$ and $\varphi_{2} \not \equiv 0$. According to (\ref{equ:1.2}), we
pick $ \tau \in \textrm{Rad}_{\phi} \bigcap (\bigcap _{2 \not= p \in
\overline{1,2n}} \textrm{ker}_{\varphi_{p}}) \backslash
\textrm{ker}_{\varphi_{2}}$ and  $ \tau' \in \textrm{Rad}_{\phi}
\bigcap (\bigcap _{ p \in \overline{2,2n}} \textrm{ker}_{\varphi_p})
\backslash \textrm{ker}_{\varphi_{1}}$ such that
\begin{equation}
\varphi_{1}(\alpha+\tau')+1 \not=0 \textrm{ and }
\varphi_{1}(\alpha+\tau')+2 \not=0.
\end{equation}
Then we find
\begin{eqnarray}
\mathcal{I}_{[0]}\cap \mathscr{A}_{\langle \hat{k} \rangle} &\owns&
[x^{-\tau, \mathbf{0}},[x^{\tau+\tau', \mathbf{0}}, u]]\equiv
\sum_{\substack{(-\varepsilon, \mathbf{0})\not =(\gamma,
\mathbf{l})\in
\textrm{ker}_{\varphi_{2}} \times \boldsymbol{\mathcal{J}}\\
l_{2}=0, |\mathbf{l}|= \hat{k}}} -a_{\gamma,
\mathbf{l}}\varphi_{2}^{2}(\tau)
(\varphi_{1}(\gamma)+1)\cdot\nonumber\\
 & &(\varphi_{1}(\gamma+\tau')+2)
x^{\gamma+\tau', \mathbf{l}} \;\;(\textrm{mod}\ \mathscr{A}_{\langle
\hat{k}-1 \rangle}).\label{equ:2.22}
\end{eqnarray}
Since $\iota(u)$ is minimal, $\alpha+\tau'\not=-\varepsilon$ due to
$\varphi_{1}(\alpha+ \tau')+1 \not=0$, and
\begin{equation}
a_{\alpha, \mathbf{i}} \varphi_{2}^{2}(\tau)
(\varphi_{1}(\alpha)+1)(\varphi_{1}(\alpha+\tau')+2)\not=0,
\end{equation}
we have
\begin{equation}
a_{\beta, \mathbf{j}} \varphi_{2}^{2}(\tau)
(\varphi_{1}(\beta)+1)(\varphi_{1}(\beta+\tau')+2)\not=0.
\end{equation}
But
\begin{equation}\label{equ:2.23}
\varphi_{1}(\alpha+\tau')\not=\varphi_{1}(\beta+\tau'),
\end{equation}
and
\begin{equation}\label{equ:2.24}
\varphi_{1}(\alpha+\tau'+\beta+\tau')+3 = 2 \varphi_{1}(\tau') \not=
0,
\end{equation}
which contradicts (\ref{equ:2.21}) if we replace $u$ by $[x^{-\tau,
\mathbf{0}},[x^{\tau+\tau', \mathbf{0}}, u]]$. Thus the first
equation in (\ref{equ:2.21}) holds, and so does (\ref{equ:2.15}).

Secondly, we want to prove that
\begin{equation}\label{equ:2.25}
\varphi_{q}(\alpha)=\varphi_{q}(\beta) \quad \textrm{whenever}\
a_{\alpha, \mathbf{i}}a_{\beta, \mathbf{j}} \not= 0 \ \textrm{for }
q\in\overline{3,2n}.
\end{equation}

Fixing $q$, if $\varphi_{q} \equiv 0$, it is done.  For
$s\in\overline{0,n-1}$, we define
\begin{equation}
(2s+1)'=2s+2,\ (2s+2)'=2s+1,\ \epsilon_{2s+1}=1,\
\epsilon_{2s+2}=-1.\end{equation} Assume $\varphi_{q} \not \equiv 0$
and $\varphi_{q'} \not \equiv 0$. Then for any  $0\neq \tau \in
\textrm{Rad}_{\phi} \cap (\cap _{q' \not= p \in \overline{1,2n}}
\textrm{ker}_{\varphi_p})$, we have
\begin{eqnarray}\label{equ:2.26}
&&\mathcal{I}_{[0]}\cap \mathscr{A}_{\langle \hat{k} \rangle} \owns
[x^{-\tau - 2\varepsilon - 2\sigma_{q}, \mathbf{0}}, [u, x^{\tau, \mathbf{0}}]]\nonumber\\
 & \equiv & \sum_{\substack{(-\varepsilon, \mathbf{0})\not =(\alpha,
\mathbf{i})\in
\textrm{ker}_{\varphi_{2}} \times \boldsymbol{\mathcal{J}}\\
i_{2}=0, |\mathbf{i}|= \hat{k}}} a_{\alpha, \mathbf{i}}
\varphi_{q'}(\tau) \varphi_{q}(\alpha)
(\varphi_{q}(\alpha+\sigma_{q}) \varphi_{q'}(\tau + 2\sigma_{q}) \nonumber \\
& & - 2 \varphi_{q} (\sigma_{q}) \varphi_{q'} (\alpha + \tau +
\sigma_{q}))x^{\alpha, \mathbf{i}} \;\;(\textrm{mod}\
\mathscr{A}_{\langle \hat{k}-1 \rangle})
\end{eqnarray}
(cf. (1.10)). Since $\varphi_{q'}(\tau)$ takes an infinite number of
elements in $\mathbb{F}$ because $\mathbb{Z}\tau\subset
\textrm{Rad}_{\phi} \cap (\cap _{q' \not= p \in \overline{1,2n}}
\textrm{ker}_{\varphi_p})$, the coefficients of
$\varphi_{q'}^{2}(\tau)$ in (\ref{equ:2.26}) show
\begin{equation}\label{equ:2.27}
\varphi_{q}(\alpha)\varphi_{q}(\alpha-\sigma_{q})=\varphi_{q}(\beta)\varphi_{q}(\beta-\sigma_{q})
\qquad \textrm{whenever}\ a_{\alpha, \mathbf{i}}a_{\beta,
\mathbf{j}} \not= 0
\end{equation}
by the minimality of $\iota (u)$ and Lemma \ref{lemma:1}. Moreover,
(\ref{equ:2.27}) is equivalent to
\begin{equation}\label{equ:2.28}
\varphi_{q}(\alpha)=\varphi_{q}(\beta) \ \textrm{or}\
\varphi_{q}(\alpha+\beta-\sigma_{q})=0.
\end{equation}
Assume that there exist $\alpha, \beta \in \Gamma$ such that
$\varphi_{q}(\alpha) \not= \varphi_{q}(\beta),\
\varphi_{q}(\alpha+\beta-\sigma_{q})=0$ and $a_{\alpha,
\mathbf{i}}a_{\beta, \mathbf{j}} \not= 0$. Without loss of
generality, we assume that $\varphi_{q}(\alpha) \not= 0$. Since
$\varphi_{q} \not \equiv 0$ and $\varphi_{q'} \not \equiv 0$, we can
choose
\begin{equation}
\tau \in \textrm{Rad}_{\phi} \bigcap (\bigcap _{q \not= p \in
\overline{1,2n}} \textrm{ker}_{\varphi_p}) \backslash
\textrm{ker}_{\varphi_q},\;\; \tau' \in \textrm{Rad}_{\phi} \bigcap
(\bigcap _{q' \not= p \in \overline{1,2n}} \textrm{ker}_{\varphi_p})
\backslash \textrm{ker}_{\varphi_{q'}}
\end{equation}
by (\ref{equ:1.2}) such that
\begin{equation}\label{equ:2.29}
\varphi_{q}(\tau+\alpha) \not= 0, \quad
\varphi_{q}(\alpha)\varphi_{q'}(\tau' -\sigma_{q}) \not=
\varphi_{q'}(\alpha)\varphi_{q}(\tau -\sigma_{q}).
\end{equation}
Then we have
\begin{eqnarray}
\mathcal{I}_{[0]}\cap \mathscr{A}_{\langle \hat{k} \rangle} &\owns&
[x^{\tau+\tau' -\varepsilon-\sigma_{q}, \mathbf{0}}, u]
  \equiv  \sum_{\substack{(-\varepsilon, \mathbf{0})\not =(\gamma,
\mathbf{l})\in
\textrm{ker}_{\varphi_{2}} \times \boldsymbol{\mathcal{J}}\\
l_{2}=0, |\mathbf{l}|= \hat{k}}} a_{\gamma, \mathbf{l}} \epsilon_{q}
(\varphi_{q}(\tau-\sigma_{q})\varphi_{q'}(\gamma) \nonumber\\&
&-\varphi_{q'}(\tau'-\sigma_{q}) \varphi_{q}(\gamma))
x^{\gamma+\tau+\tau', \mathbf{l}} \;\;(\textrm{mod}\
\mathscr{A}_{\langle \hat{k}-1 \rangle}).\label{equ:2.30}
\end{eqnarray}
Since $\iota(u)$ is minimal, $\alpha+\tau+\tau'\not=-\varepsilon$
due to $\varphi_{q}(\alpha +\tau+ \tau') =\varphi_{q}(\alpha +\tau)
\not=0$, and
\begin{equation}
a_{\alpha, \mathbf{i}} \epsilon_{q}
(\varphi_{q}(\tau-\sigma_{q})\varphi_{q'}(\alpha)
-\varphi_{q'}(\tau'-\sigma_{q}) \varphi_{q}(\alpha))\not=0,
\end{equation}
we have
\begin{equation}
a_{\beta, \mathbf{j}} \epsilon_{q}
(\varphi_{q}(\tau-\sigma_{q})\varphi_{q'}(\beta)
-\varphi_{q'}(\tau'-\sigma_{q}) \varphi_{q}(\beta))\not=0.
\end{equation}
But
\begin{equation}\label{equ:2.31}
\varphi_{q}(\alpha+\tau+\tau')\not=\varphi_{q}(\beta+\tau+\tau') \
\textrm{and}\
\varphi_{q}(\alpha+\tau+\tau'+\beta+\tau+\tau'-\sigma_{q})\not= 0,
\end{equation}
which contradicts (\ref{equ:2.28}) if we replace $u$ by
$[x^{\tau+\tau' -\varepsilon - \sigma_{q}, \mathbf{0}}, u]$. Thus
the first equation in (\ref{equ:2.28}) holds, and so does
(\ref{equ:2.25}).

If $\varphi_{q} \not \equiv 0$ and $\varphi_{q'} \equiv 0$, then
$\mathcal {J}_{q'} = \mathbb{N}$ by (1.7). Moreover, if $j_{q'}\not
= 0$ for some $a_{\beta, \mathbf{j}} \not = 0$, choose $ \tau \in
\textrm{Rad}_{\phi}\cap(\cap _{q \not= p \in \overline{1,2n}}
\textrm{ker}_{\varphi_p} )\backslash \textrm{ker}_{\varphi_{q}}$
such that $\varphi_{q}(\beta+\tau+\sigma_{q}) \not = 0.$ Then we
have
\begin{eqnarray}
&&\mathcal{I}_{[0]}\cap \mathscr{A}_{\langle \hat{k}-1 \rangle}
\owns
[x^{\tau, \mathbf{0}}, u]\nonumber\\
 & \equiv & \sum_{\substack{(-\varepsilon, \mathbf{0})\not =(\gamma,
\mathbf{l})\in
\textrm{ker}_{\varphi_{2}} \times \boldsymbol{\mathcal{J}}\\
l_{2}=0, |\mathbf{l}|= \hat{k}}} a_{\gamma, \mathbf{l}} \epsilon_{q}
l_{q'}\varphi_{q}(\tau) x^{\gamma+\tau+\sigma_{q}+\varepsilon,
\mathbf{l}-1_{[q']}} \;\;(\textrm{mod}\ \mathscr{A}_{\langle
\hat{k}-2 \rangle}).\label{equ:2.32}
\end{eqnarray}
Since $\varphi_{q}(\beta+\tau+\sigma_{q}+\varepsilon) =
\varphi_{q}(\beta+\tau+\sigma_{q}) \not = 0$ and $j_{q'}\not = 0$
for some $a_{\beta, \mathbf{j}} \not = 0$, we have $[x^{\tau,
\mathbf{0}}, u] \in (\mathcal{I}_{[0]}\cap \mathscr{A}_{\langle
\hat{k}-1 \rangle}) \backslash \mathbb{F} x^{-\varepsilon,
\mathbf{0}}$, which contradicts the minimality of $\hat{k}$.
Therefore $j_{q'} = 0$ whenever $a_{\beta, \mathbf{j}} \not = 0$.
Thus we get
\begin{eqnarray}
&&\mathcal{I}_{[0]}\cap \mathscr{A}_{\langle \hat{k} \rangle} \owns
[x^{-\sigma_{q}-\varepsilon, 1_{[q']}}, u]\nonumber\\
 & \equiv & \sum_{\substack{(-\varepsilon, \mathbf{0})\not =(\gamma,
\mathbf{l})\in
\textrm{ker}_{\varphi_{2}} \times \boldsymbol{\mathcal{J}}\\
l_{2}=0, |\mathbf{l}|= \hat{k}}} -a_{\gamma, \mathbf{l}}
\epsilon_{q} \varphi_{q}(\gamma) x^{\gamma, \mathbf{l}}
\;\;(\textrm{mod}\ \mathscr{A}_{\langle \hat{k}-1
\rangle}),\label{equ:2.33}
\end{eqnarray}
which implies (\ref{equ:2.25}). Therefore, we have proved
\begin{equation}\label{equ:2.34}
\varphi_{q}(\alpha)=\varphi_{q}(\beta) \quad \textrm{whenever}\
a_{\alpha, \mathbf{i}}a_{\beta, \mathbf{j}} \not= 0, \quad
\textrm{for } q \in \overline{1,2n}.
\end{equation}

Thirdly, we want to prove that
\begin{equation}\label{equ:2.35}
\phi(\gamma, \alpha)=\phi(\gamma, \beta)\textrm{ for any }\gamma \in
\Gamma,\textrm{ whenever } a_{\alpha, \mathbf{i}}a_{\beta,
\mathbf{j}}\not=0.
\end{equation}

For $a_{\beta, \mathbf{j}} \not = 0$, we have
\begin{eqnarray}
&&\mathcal{I}_{[0]}\cap \mathscr{A}_{\langle \hat{k} \rangle} \owns
[x^{\beta, \mathbf{0}}, u]\nonumber\\
 & \equiv & \sum_{\substack{(-\varepsilon, \mathbf{0})\not =(\alpha,
\mathbf{i})\in
\textrm{ker}_{\varphi_{2}} \times \boldsymbol{\mathcal{J}}\\
i_{2}=0, |\mathbf{i}|= \hat{k}}} a_{\alpha, \mathbf{i}} \phi (\beta,
\alpha) x^{\alpha+\beta+\varepsilon, \mathbf{i}}\;\; (\textrm{mod}\
\mathscr{A}_{\langle \hat{k}-1 \rangle})\label{equ:4.1}
\end{eqnarray}
by (\ref{equ:2.1}) and (\ref{equ:2.34}). Since $\iota(u)$ is minimal
and $\phi(\beta, \beta)=0$, we have
\begin{equation}\label{equ:4.2}
\phi(\alpha, \beta)=0 \qquad \textrm{whenever}\ a_{\alpha,
\mathbf{i}}a_{\beta, \mathbf{j}} \not= 0.
\end{equation}

For any $\gamma \in \Gamma$, we have
\begin{eqnarray}
&&\mathcal{I}_{[0]}\cap \mathscr{A}_{\langle \hat{k} \rangle} \owns
[x^{-\gamma-\varepsilon,\mathbf{0}},[x^{\gamma-\varepsilon, \mathbf{0}}, u]]\nonumber\\
 & \equiv & \sum_{\substack{(-\varepsilon, \mathbf{0})\not =(\alpha,
\mathbf{i})\in
\textrm{ker}_{\varphi_{2}} \times \boldsymbol{\mathcal{J}}\\
i_{2}=0, |\mathbf{i}|= \hat{k}}} a_{\alpha, \mathbf{i}} \Big(
\sum_{s=1}^{n-1} \sum_{t=1}^{n-1}
(\varphi_{2s+1}(\gamma)\varphi_{2s+2}(\alpha) -
\varphi_{2s+1}(\alpha)\varphi_{2s+2}(\gamma)) \cdot \nonumber\\
& & (\varphi_{2t+1}(\alpha+\sigma_{2s+1})\varphi_{2t+2}(\gamma) -
\varphi_{2t+1}(\gamma)\varphi_{2t+2}(\alpha+\sigma_{2s+1}))
x^{\alpha+\sigma_{2s+1}+\sigma_{2t+1}, \mathbf{i}} \nonumber \\
& & +2 \sum_{s=1}^{n-1}
(\varphi_{2s+1}(\gamma)\varphi_{2s+2}(\alpha) -
\varphi_{2s+1}(\alpha)\varphi_{2s+2}(\gamma))
\varphi_{2}(\gamma) (\varphi_{1}(\alpha)+1) x^{\alpha+\sigma_{2s+1}-\varepsilon, \mathbf{i}} \nonumber \\
& & -2 \sum_{s=1}^{n-1}
(\varphi_{2s+1}(\gamma)\varphi_{2s+2}(\alpha) -
\varphi_{2s+1}(\alpha)\varphi_{2s+2}(\gamma)) \phi(\gamma, \alpha) x^{\alpha+\sigma_{2s+1}, \mathbf{i}} \nonumber \\
& & -\phi^{2}(\gamma,\alpha)x^{\alpha, \mathbf{i}}  +2\phi(\gamma,\alpha)\varphi_{2}(\gamma) (\varphi_{1}(\alpha)+1)x^{\alpha-\varepsilon, \mathbf{i}} \nonumber \\
& & - \varphi_{2}^{2}(\gamma) \varphi_{1}(\alpha)
(\varphi_{1}(\alpha)+1)x^{\alpha-2\varepsilon, \mathbf{i}} \Big)\ \
(\textrm{mod}\ \mathscr{A}_{\langle \hat{k}-1
\rangle}).\label{equ:2.36}
\end{eqnarray}

For $a_{\alpha, \mathbf{i}} \not = 0 $, if $ \varphi_{q}(\alpha)=0$
for $q \in\overline{3,2n}$ and $ \varphi_{1}(\alpha)+1=0$, then
(\ref{equ:2.34}) and (\ref{equ:2.36}) imply
\begin{eqnarray}
&&\mathcal{I}_{[0]}\cap \mathscr{A}_{\langle \hat{k} \rangle} \owns
[x^{-\gamma-\varepsilon,\mathbf{0}},[x^{\gamma-\varepsilon, \mathbf{0}}, u]]\nonumber\\
 & \equiv & \sum_{\substack{(-\varepsilon, \mathbf{0})\not =(\alpha,
\mathbf{i})\in
\textrm{ker}_{\varphi_{2}} \times \boldsymbol{\mathcal{J}}\\
i_{2}=0, |\mathbf{i}|= \hat{k}}} -a_{\alpha, \mathbf{i}}
\phi^{2}(\gamma,\alpha)x^{\alpha, \mathbf{i}} \quad(\textrm{mod}\
\mathscr{A}_{\langle \hat{k}-1 \rangle}).\label{equ:2.37}
\end{eqnarray}
By the minimality of $\iota (u)$ and Lemma \ref{lemma:1}, we have
\begin{equation}\label{equ:2.38}
\phi^{2}(\gamma,\alpha) = \phi^{2}(\gamma,\beta) \quad
\textrm{whenever} \ a_{\alpha, \mathbf{i}}a_{\beta, \mathbf{j}}
\not= 0,
\end{equation}
which implies
\begin{equation}\label{equ:2.39}
\phi(\gamma,\alpha) = \phi(\gamma,\beta)\textrm{ or }
\phi(\gamma,\alpha) = -\phi(\gamma,\beta) \ \quad \textrm{whenever}\
a_{\alpha, \mathbf{i}}a_{\beta, \mathbf{j}} \not= 0.
\end{equation}
Assume that there exist $\alpha, \beta \in \Gamma$ such that
\begin{equation}
\phi(\gamma,\alpha) \not = \phi(\gamma,\beta),\;\;
\phi(\gamma,\alpha) = -\phi(\gamma,\beta) \textrm{ and } a_{\alpha,
\mathbf{i}}a_{\beta, \mathbf{j}} \not= 0.
\end{equation}
Then we
obtain $\phi(\gamma,\alpha) = -\phi(\gamma,\beta) \not= 0$.
Moreover,  (\ref{equ:2.34}) and (\ref{equ:4.2}) tell us that
\begin{eqnarray}
&&\mathcal{I}_{[0]}\cap \mathscr{A}_{\langle \hat{k} \rangle} \owns
[x^{-\gamma+\alpha,\mathbf{0}},[x^{\gamma-\varepsilon, \mathbf{0}},
u]]\nonumber\\ &\equiv &- \sum_{\substack{(-\varepsilon,
\mathbf{0})\not =(\theta, \mathbf{i})\in
\textrm{ker}_{\varphi_{2}} \times \boldsymbol{\mathcal{J}}\\
i_{2}=0, |\mathbf{i}|= \hat{k}}} a_{\theta, \mathbf{i}}
\phi(\gamma,\theta)
\phi(\gamma,\theta+\alpha)x^{\theta+\alpha+\varepsilon, \mathbf{i}}\
(\textrm{mod}\ \mathscr{A}_{\langle \hat{k}-1
\rangle}).\label{equ:2.40}
\end{eqnarray}
Since $\phi(\gamma,2\alpha+\varepsilon)=2\phi(\gamma,\alpha)\not=0$,
we have $2\alpha+\varepsilon\not=-\varepsilon$. But the fact that
\begin{equation}
a_{\alpha, \mathbf{i}}
\phi(\gamma,\alpha)\phi(\gamma,2\alpha)\not=0,\;\; a_{\beta,
\mathbf{j}} \phi(\gamma,\beta)\phi(\gamma,\alpha+\beta)=0
\end{equation}
contradicts the minimality of $\iota(u)$. Thus the first equation in
(\ref{equ:2.39}) holds, and so does (\ref{equ:2.35}).

Next we assume that
\begin{equation}\label{equ:2.41}
 \varphi_{1}(\alpha) + 1\not=0 \textrm{ or there exists } q \in\overline{3,2n} \textrm{ such that } \varphi_{q}(\alpha)
\not=0, \textrm{ for } a_{\alpha, \mathbf{i}} \not = 0 .
\end{equation}

\noindent{\emph{Claim 1.}} There exists some $v_{1} \in
\mathcal{I}_{[0]}\cap \mathscr{A}_{\langle \hat{k} \rangle} $ such
that
\begin{eqnarray}
v_{1} & \equiv & \sum_{\substack{(-\varepsilon, \mathbf{0})\not
=(\alpha, \mathbf{i})\in
\textrm{ker}_{\varphi_{2}} \times \boldsymbol{\mathcal{J}}\\
i_{2}=0, |\mathbf{i}|= \hat{k}}} a_{\alpha, \mathbf{i}} \Big(
\sum_{s=1}^{n-1} \sum_{t=1}^{n-1}
(\varphi_{2s+1}(\gamma)\varphi_{2s+2}(\alpha) -
\varphi_{2s+1}(\alpha)\varphi_{2s+2}(\gamma)) \cdot \nonumber\\
& & (\varphi_{2t+1}(\alpha+\sigma_{2s+1})\varphi_{2t+2}(\gamma) -
\varphi_{2t+1}(\gamma)\varphi_{2t+2}(\alpha+\sigma_{2s+1}))
x^{\alpha+\sigma_{2s+1}+\sigma_{2t+1}, \mathbf{i}} \nonumber \\
& & -2 \sum_{s=1}^{n-1}
(\varphi_{2s+1}(\gamma)\varphi_{2s+2}(\alpha) -
\varphi_{2s+1}(\alpha)\varphi_{2s+2}(\gamma)) \phi(\gamma, \alpha) x^{\alpha+\sigma_{2s+1}, \mathbf{i}} \nonumber \\
& & -\phi^{2}(\gamma,\alpha)x^{\alpha, \mathbf{i}} \Big)\;\;
(\textrm{mod}\ \mathscr{A}_{\langle \hat{k}-1
\rangle}).\label{equ:2.43}
\end{eqnarray}

For convenience, we denote
\begin{eqnarray}
U' & := & \mathscr{A}_{\langle \hat{k} \rangle}\cap (\sum_{(\alpha,
\mathbf{i})\in \textrm{ker}_{\varphi_{2}} \times
\boldsymbol{\mathcal{J}}; i_{2}=0} \mathbb{F} x^{\alpha,
\mathbf{i}}),\\
W_{1} & := & \mathscr{A}_{\langle \hat{k}-1 \rangle}\cap
(\sum_{(\alpha, \mathbf{i})\in \textrm{ker}_{\varphi_{2}} \times
\boldsymbol{\mathcal{J}}; i_{2}=0} \mathbb{F} x^{\alpha,
\mathbf{i}}),\\
W_{2} & := & \mathcal{I}_{[0]}\cap \mathscr{A}_{\langle \hat{k}
\rangle}=\mathcal{I}\cap U'.
\end{eqnarray}
Notice that $W_{1}$ and $W_{2}$ are two subspaces of $U'$. If
$\mathcal{J}_{2} = \mathbb{N}$, $\textrm{ad}_{x^{0,1_{[2]}}}$ is a
linear transformation on $\mathscr{A}$ which preserves $U'$, $W_{1}$
and $W_{2}$. For  any $a_{\alpha, \mathbf{i}} \not = 0 $, we have
\begin{equation}\label{equ:2.42}
[x^{0,1_{[2]}}, x^{\alpha, \mathbf{i}}] \equiv -
(\varphi_{1}(\alpha)+1) x^{\alpha, \mathbf{i}}\;\;(\textrm{mod}\
W_{1}).
\end{equation}
Thus, acting  on $U'/W_{1}$, $\textrm{ad}_{x^{0,1_{[2]}}}$ has
eigenvectors
$x^{\alpha+\sigma_{2s+1}+\sigma_{2t+1},\mathbf{i}}+W_{1}$,
$x^{\alpha+\sigma_{2s+1}, \mathbf{i}}+W_{1}$ and $x^{\alpha,
\mathbf{i}}+W_{1}$ with eigenvalue $- \varphi_{1}(\alpha)-1$.
Moreover, it has eigenvectors $x^{\alpha+\sigma_{2s+1}-\varepsilon,
\mathbf{i}}+W_{1}$ and $x^{\alpha-\varepsilon, \mathbf{i}}+W_{1}$
with eigenvalue $- \varphi_{1}(\alpha)$. Furthermore,
$x^{\alpha-2\varepsilon, \mathbf{i}}+W_{1}$ is an eigenvector with
eigenvalue $- \varphi_{1}(\alpha)+1$. Repeatedly applying
$\textrm{ad}_{x^{0,1_{[2]}}}$ on $(W_{2}+W_{1})/W_{1}$, \emph{Claim
1} holds by Lemma \ref{lemma:1}, (\ref{equ:2.34}) and
(\ref{equ:2.36}).

If $\mathcal{J}_{2}=\{0\}$, we have $\varphi_{2}\not\equiv 0$ by
(1.7). For any $a_{\alpha, \mathbf{i}} \not = 0 $, if
$\varphi_{1}(\alpha)+1=0$, \emph{Claim 1} follows naturally from
(\ref{equ:2.36}). Suppose that $\varphi_{1}(\alpha)+1\not=0$.
Picking $\tau \in \textrm{Rad}_{\phi} \cap (\cap_{2\not=p \in
\overline{1,2n}} \textrm{ker}_{\varphi_{p}}) \backslash
\textrm{ker}_{\varphi_{2}}$, $\textrm{ad}_{x^{-\tau,
\mathbf{0}}}\textrm{ad}_{x^{\tau, \mathbf{0}}}$ is a linear
transformation on $\mathscr{A}$ which preserves $U'$, $W_{1}$ and
$W_{2}$. Since for any $a_{\alpha, \mathbf{i}} \not = 0 $, we have
\begin{equation}\label{equ:2.45}
 [x^{-\tau, \mathbf{0}},[x^{\tau, \mathbf{0}}, x^{\alpha,
\mathbf{i}}]]
  \equiv g(\alpha)
  x^{\alpha,
\mathbf{i}}\  (\textrm{mod}\ W_{1}),
\end{equation}
where $g(\alpha)=
-\varphi^{2}_{2}(\tau)(\varphi_{1}(\alpha)+1)(\varphi_{1}(\alpha)+2)$.
Acting on $U'/W_{1}$, $\textrm{ad}_{x^{-\tau,
\mathbf{0}}}\textrm{ad}_{x^{\tau, \mathbf{0}}}$  has eigenvectors
$x^{\alpha+\sigma_{2s+1}+\sigma_{2t+1}, \mathbf{i}}+W_{1},\
x^{\alpha+\sigma_{2s+1}, \mathbf{i}}+W_{1}$ and $x^{\alpha,
\mathbf{i}}+W_{1}$ with eigenvalue $g(\alpha)$. Moreover, it has
eigenvectors $x^{\alpha+\sigma_{2s+1}-\varepsilon,
\mathbf{i}}+W_{1}$ and $ x^{\alpha-\varepsilon, \mathbf{i}}+W_{1}$
with eigenvalue $g(\alpha-\varepsilon)$. Furthermore,
$x^{\alpha-2\varepsilon, \mathbf{i}}+W_{1}$ is an eigenvector with
eigenvalue $g(\alpha-2\varepsilon)$. First, we have
$g(\alpha)\not=g(\alpha-\varepsilon)$. If $2\varphi_{1}(\alpha)+1
\not=0$, we further get $g(\alpha)\not=g(\alpha-2\varepsilon)$.
Thus, repeatedly applying $\textrm{ad}_{x^{-\tau,
\mathbf{0}}}\textrm{ad}_{x^{\tau, \mathbf{0}}}$ on
$(W_{2}+W_{1})/W_{1}$, \emph{Claim 1} holds by Lemma \ref{lemma:1},
(\ref{equ:2.34}) and (\ref{equ:2.36}). Otherwise,
$2\varphi_{1}(\alpha)+1 =0$, we have
$g(\alpha)=g(\alpha-2\varepsilon)$. By Lemma \ref{lemma:1},
(\ref{equ:2.34}) and (\ref{equ:2.36}) we obtain some $v_{1}' \in
\mathcal{I}_{[0]}\cap \mathscr{A}_{\langle \hat{k} \rangle} $ that
\begin{eqnarray}
v_{1}' & \equiv & \sum_{\substack{(-\varepsilon, \mathbf{0})\not
=(\alpha, \mathbf{i})\in
\textrm{ker}_{\varphi_{2}} \times \boldsymbol{\mathcal{J}}\\
i_{2}=0, |\mathbf{i}|= \hat{k}}} a_{\alpha, \mathbf{i}} \Big(
\sum_{s=1}^{n-1} \sum_{t=1}^{n-1}
(\varphi_{2s+1}(\gamma)\varphi_{2s+2}(\alpha) -
\varphi_{2s+1}(\alpha)\varphi_{2s+2}(\gamma)) \cdot \nonumber\\
& & (\varphi_{2t+1}(\alpha+\sigma_{2s+1})\varphi_{2t+2}(\gamma) -
\varphi_{2t+1}(\gamma)\varphi_{2t+2}(\alpha+\sigma_{2s+1}))
x^{\alpha+\sigma_{2s+1}+\sigma_{2t+1}, \mathbf{i}} \nonumber \\
& & -2 \sum_{s=1}^{n-1}
(\varphi_{2s+1}(\gamma)\varphi_{2s+2}(\alpha) -
\varphi_{2s+1}(\alpha)\varphi_{2s+2}(\gamma)) \phi(\gamma, \alpha) x^{\alpha+\sigma_{2s+1}, \mathbf{i}} \nonumber \\
& & -\phi^{2}(\gamma,\alpha)x^{\alpha, \mathbf{i}}-
\varphi_{2}^{2}(\gamma) \varphi_{1}(\alpha)
(\varphi_{1}(\alpha)+1)x^{\alpha-2\varepsilon, \mathbf{i}}
\Big)\quad(\textrm{mod}\ \mathscr{A}_{\langle \hat{k}-1
\rangle}).\label{equ:2.46}
\end{eqnarray}
Since $2\varphi_{1}(\alpha)+1 =0$, we have $\varphi_{1}(\alpha)\not
=0$. Pick $\tau \in \textrm{Rad}_{\phi} \cap (\cap_{2\not=p \in
\overline{1,2n}} \textrm{ker}_{\varphi_{p}}) \backslash
\textrm{ker}_{\varphi_{2}}$ again, then we have
\begin{eqnarray}
\mathcal{I}_{[0]}\cap \mathscr{A}_{\langle \hat{k} \rangle} &\owns&
[x^{-\tau, \mathbf{0}},[x^{\tau-2\varepsilon, \mathbf{0}}, u]]
  \equiv  \sum_{\substack{(-\varepsilon, \mathbf{0})\not =(\alpha,
\mathbf{i})\in
\textrm{ker}_{\varphi_{2}} \times \boldsymbol{\mathcal{J}}\\
i_{2}=0, |\mathbf{i}|= \hat{k}}} -a_{\alpha, \mathbf{i}}
\varphi^{2}_{2}(\tau)
\varphi_{1}(\alpha)\cdot\nonumber\\
&& (\varphi_{1}(\alpha)+1) x^{\alpha-2\varepsilon, \mathbf{i}}\
(\textrm{mod}\ \mathscr{A}_{\langle \hat{k}-1
\rangle}).\label{equ:2.47}
\end{eqnarray}
Therefore, \emph{Claim 1} holds by substituting (\ref{equ:2.46}) and
(\ref{equ:2.47}) in
\begin{equation}
v_{1}'-\frac{\varphi_{2}^{2}(\gamma)}{\varphi^{2}_{2}(\tau)}[x^{-\tau,
\mathbf{0}},[x^{\tau-2\varepsilon, \mathbf{0}}, u]].
\end{equation}

\noindent{\emph{Claim 2.}} For any fixed $r$ (where $r \in
\overline{1,n-1}$), there exists some $v_{2} \in
\mathcal{I}_{[0]}\cap \mathscr{A}_{\langle \hat{k} \rangle} $ such
that
\begin{eqnarray}
v_{2} & \equiv & \sum_{\substack{(-\varepsilon, \mathbf{0})\not
=(\alpha, \mathbf{i})\in
\textrm{ker}_{\varphi_{2}} \times \boldsymbol{\mathcal{J}}\\
i_{2}=0, |\mathbf{i}|= \hat{k}}} a_{\alpha, \mathbf{i}} \Big(
\sum_{\substack{s=1\\s\not=r}}^{n-1}
\sum_{\substack{t=1\\t\not=r}}^{n-1}
(\varphi_{2s+1}(\gamma)\varphi_{2s+2}(\alpha) -
\varphi_{2s+1}(\alpha)\varphi_{2s+2}(\gamma)) \cdot \nonumber\\
& & (\varphi_{2t+1}(\alpha+\sigma_{2s+1})\varphi_{2t+2}(\gamma) -
\varphi_{2t+1}(\gamma)\varphi_{2t+2}(\alpha+\sigma_{2s+1}))
x^{\alpha+\sigma_{2s+1}+\sigma_{2t+1}, \mathbf{i}} \nonumber \\
& & -2 \sum_{\substack{s=1\\s\not=r}}^{n-1}
(\varphi_{2s+1}(\gamma)\varphi_{2s+2}(\alpha) -
\varphi_{2s+1}(\alpha)\varphi_{2s+2}(\gamma)) \phi(\gamma, \alpha) x^{\alpha+\sigma_{2s+1}, \mathbf{i}} \nonumber \\
& & -\phi^{2}(\gamma,\alpha)x^{\alpha, \mathbf{i}}
\Big)\quad(\textrm{mod}\ \mathscr{A}_{\langle \hat{k}-1
\rangle}).\label{equ:2.49}
\end{eqnarray}

If $\sigma_{2r+1}=0$, we get
$\varphi_{2r+1}\equiv\varphi_{2r+2}\equiv0$ by (\ref{equ:1.4}).
Therefore \emph{Claim 2} follows naturally from (\ref{equ:2.43}). We
assume that $\sigma_{2r+1} \not= 0$. Then
$\varphi_{2r+1}(\sigma_{2r+1}) \not= 0$ or
$\varphi_{2r+2}(\sigma_{2r+1}) \not= 0$ by (\ref{equ:1.3}) and
(\ref{equ:1.4}). Choose $q \in \{2r+1,\,2r+2\}$ such that
$\varphi_{q}(\sigma_{2r+1}) \not= 0$, then $\varphi_{q} \not\equiv
0$ and $\sigma_{q}=\sigma_{2r+1}$.

If $\varphi_{q'} \equiv 0$, we get $\mathcal {J}_{q'} = \mathbb{N}$
by (\ref{equ:1.1}). Observe that
$\textrm{ad}_{x^{-\sigma_{q}-\varepsilon, 1_{[q']}}}$ is a linear
transformation on $\mathscr{A}$ which preserves $U'$, $W_{1}$ and
$W_{2}$. For any $a_{\alpha, \mathbf{i}} \not = 0 $, since $i_{q'} =
0$ by (\ref{equ:2.32}), we have
\begin{equation}\label{equ:2.48}
[x^{-\sigma_{q}-\varepsilon, 1_{[q']}}, x^{\alpha, \mathbf{i}}]
  \equiv   - \epsilon_{q} \varphi_{q}(\alpha)
x^{\alpha, \mathbf{i}} \  (\textrm{mod}\ W_{1}).
\end{equation}
Repeatedly applying $\textrm{ad}_{x^{-\sigma_{q}-\varepsilon,
1_{[q']}}}$ on $(W_{2}+W_{1})/W_{1}$, \emph{Claim 2} holds by Lemma
\ref{lemma:1}, (\ref{equ:2.34}), (\ref{equ:2.43}) and
(\ref{equ:2.48}).

If $\varphi_{q'} \not\equiv 0$, pick $ \tau \in \textrm{Rad}_{\phi}
\cap (\cap _{q' \not= p \in \overline{1,2n}}
\textrm{ker}_{\varphi_p})\backslash \textrm{ker}_{\varphi_{q'}}$.
Then $\textrm{ad}_{x^{-\tau - 2\varepsilon - 2\sigma_{q},
\mathbf{0}}}\textrm{ad}_{x^{\tau, \mathbf{0}}}$ is a linear
transformation on $\mathscr{A}$ which preserves $U'$, $W_{1}$ and
$W_{2}$. For any $a_{\alpha, \mathbf{i}} \not = 0 $, we have
\begin{equation}\label{equ:2.50}
[x^{-\tau - 2\varepsilon - 2\sigma_{q}, \mathbf{0}}, [x^{\tau,
\mathbf{0}}, x^{\alpha, \mathbf{i}}]] \equiv h(\alpha) x^{\alpha,
\mathbf{i}} \  (\textrm{mod}\ W_{1}),
\end{equation}
where $h(\alpha)= - \varphi_{q'}(\tau) \varphi_{q}(\alpha)
(\varphi_{q}(\alpha)\varphi_{q'}(\tau+2\sigma_{q}) -
\varphi_{q}(\sigma_{q})\varphi_{q'}(\tau+2\alpha))$. Acting on
$U'/W_{1}$, $\textrm{ad}_{x^{-\tau - 2\varepsilon - 2\sigma_{q},
\mathbf{0}}}\textrm{ad}_{x^{\tau, \mathbf{0}}}$ has eigenvectors
$x^{\alpha, \mathbf{i}}+W_{1}$,
$x^{\alpha+\sigma_{2s+1}+\sigma_{2t+1}, \mathbf{i}}+W_{1}$ and
$x^{\alpha+\sigma_{2s+1}, \mathbf{i}}+W_{1}$ (for $s, t \not=r$)
with eigenvalue $h(\alpha)$. Moreover, it has eigenvectors
$x^{\alpha+\sigma_{2r+1}, \mathbf{i}}+W_{1}$ and
$x^{\alpha+\sigma_{2s+1}+\sigma_{2r+1}, \mathbf{i}}+W_{1}$ (for
$s\not=r$) with eigenvalue $h(\alpha+\sigma_{q})$. Furthermore,
$x^{\alpha+2\sigma_{2r+1},
 \mathbf{i}}+W_{1}$ is an eigenvector with
eigenvalue $h(\alpha+2\sigma_{q})$.

For any $a_{\alpha, \mathbf{i}} \not = 0 $, if
$\varphi_{q}(\alpha)=\varphi_{q'}(\alpha)=0$, \emph{Claim 2} follows
from (\ref{equ:2.43}). Assume that $\varphi_{q}(\alpha)\not=0$ or
$\varphi_{q'}(\alpha)\not=0$. Restricting the selection of $\tau$ so
that
\begin{equation}
\varphi_{q}(\alpha)\varphi_{q'}(\tau+\sigma_{q})-\varphi_{q}(\sigma_{q})\varphi_{q'}(\alpha)\not=0,
\end{equation}
we have $h(\alpha)\not=h(\alpha+\sigma_{q})$. If
$\varphi_{q}(2\alpha+\sigma_{q})\not=0$ or
$\varphi_{q'}(2\alpha+\sigma_{q})\not=0$, we can further restrict
the selection of $\tau$ so that
\begin{equation}
\varphi_{q}(2\alpha+\sigma_{q})\varphi_{q'}(\tau+\sigma_{q})-\varphi_{q}(\sigma_{q})\varphi_{q'}(2\alpha+\sigma_{q})\not=0.
\end{equation}
Therefore we have $h(\alpha)\not=h(\alpha+2\sigma_{q})$. Repeatedly
applying $\textrm{ad}_{x^{-\tau - 2\varepsilon - 2\sigma_{q},
\mathbf{0}}}\textrm{ad}_{x^{\tau, \mathbf{0}}}$ on
$(W_{2}+W_{1})/W_{1}$, \emph{Claim 2} follows from Lemma
\ref{lemma:1}, (\ref{equ:2.34}) and (\ref{equ:2.43}). Otherwise,
$\varphi_{q}(2\alpha+\sigma_{q})=\varphi_{q'}(2\alpha+\sigma_{q})=0$,
we get $h(\alpha)=h(\alpha+2\sigma_{q})$. By Lemma \ref{lemma:1},
(\ref{equ:2.34}) and (\ref{equ:2.43}), we obtain some $v_{2}' \in
\mathcal{I}_{[0]}\cap \mathscr{A}_{\langle \hat{k} \rangle} $ such
that
\begin{eqnarray}
v_{2}' & \equiv & \sum_{\substack{(-\varepsilon, \mathbf{0})\not
=(\alpha, \mathbf{i})\in
\textrm{ker}_{\varphi_{2}} \times \boldsymbol{\mathcal{J}}\\
i_{2}=0, |\mathbf{i}|= \hat{k}}} a_{\alpha, \mathbf{i}} \Big(
\sum_{\substack{s=1\\s\not=r}}^{n-1}
\sum_{\substack{t=1\\t\not=r}}^{n-1}
(\varphi_{2s+1}(\gamma)\varphi_{2s+2}(\alpha) -
\varphi_{2s+1}(\alpha)\varphi_{2s+2}(\gamma)) \cdot \nonumber\\
& & (\varphi_{2t+1}(\alpha+\sigma_{2s+1})\varphi_{2t+2}(\gamma) -
\varphi_{2t+1}(\gamma)\varphi_{2t+2}(\alpha+\sigma_{2s+1}))
x^{\alpha+\sigma_{2s+1}+\sigma_{2t+1}, \mathbf{i}} \nonumber \\
& & -2 \sum_{\substack{s=1\\s\not=r}}^{n-1}
(\varphi_{2s+1}(\gamma)\varphi_{2s+2}(\alpha) -
\varphi_{2s+1}(\alpha)\varphi_{2s+2}(\gamma)) \phi(\gamma, \alpha) x^{\alpha+\sigma_{2s+1}, \mathbf{i}} \nonumber \\
& & + (\varphi_{2r+1}(\gamma)\varphi_{2r+2}(\alpha) -
\varphi_{2r+1}(\alpha)\varphi_{2r+2}(\gamma))^{2}
x^{\alpha+2\sigma_{2r+1}, \mathbf{i}} \nonumber\\
& & -\phi^{2}(\gamma,\alpha)x^{\alpha,
\mathbf{i}}\Big)\;\;(\textrm{mod}\ \mathscr{A}_{\langle \hat{k}-1
\rangle}).\label{equ:2.52}
\end{eqnarray}
Since $\varphi_{q}(2\alpha+\sigma_{q})=0$ and
$\varphi_{q}(\sigma_{q})\not=0$, we have $\varphi_{q}(\alpha)\not=0$
and $\varphi_{q}(\alpha+\sigma_{q})\not=0$. Picking $ \tau' \in
\textrm{Rad}_{\phi} \cap (\cap _{q' \not= p \in \overline{1,2n}}
\textrm{ker}_{\varphi_p})\backslash \textrm{ker}_{\varphi_{q'}}$, we
have
\begin{eqnarray}
\mathcal{I}_{[0]}\cap \mathscr{A}_{\langle \hat{k} \rangle} &\owns&
[x^{-\tau' - \varepsilon, \mathbf{0}}, [ x^{\tau'-\varepsilon,
\mathbf{0}},u]]
  \equiv  \sum_{\substack{(-\varepsilon, \mathbf{0})\not =(\alpha,
\mathbf{i})\in
\textrm{ker}_{\varphi_{2}} \times \boldsymbol{\mathcal{J}}\\
i_{2}=0, |\mathbf{i}|= \hat{k}}} -a_{\alpha, \mathbf{i}}
\varphi_{q'}^{2}(\tau') \varphi_{q}(\alpha)\cdot\nonumber\\
& & \varphi_{q}(\alpha+\sigma_{q})x^{\alpha +2\sigma_{q},
\mathbf{i}} \;\; (\textrm{mod}\ \mathscr{A}_{\langle \hat{k}-1
\rangle}).\label{equ:2.53}
\end{eqnarray}
Thus, \emph{Claim 2} holds by substituting (\ref{equ:2.52}) and
(\ref{equ:2.53}) in
\begin{equation}
v_{2}'+\frac{(\varphi_{2r+1}(\gamma)\varphi_{2r+2}(\alpha) -
\varphi_{2r+1}(\alpha)\varphi_{2r+2}(\gamma))^{2}}{\varphi_{q'}^{2}(\tau')
\varphi_{q}(\alpha) \varphi_{q}(\alpha+\sigma_{q})}[x^{-\tau' -
\varepsilon, \mathbf{0}}, [ x^{\tau'-\varepsilon, \mathbf{0}},u]].
\end{equation}

\noindent{\emph{Claim 3.}} There exists some $v_{n} \in
\mathcal{I}_{[0]}\cap \mathscr{A}_{\langle \hat{k} \rangle} $ such
that
\begin{equation}\label{equ:2.54}
v_{n} \equiv\sum_{\substack{(-\varepsilon, \mathbf{0})\not =(\alpha,
\mathbf{i})\in
\textrm{ker}_{\varphi_{2}} \times \boldsymbol{\mathcal{J}}\\
i_{2}=0, |\mathbf{i}|= \hat{k}}} -a_{\alpha, \mathbf{i}}
\phi^{2}(\gamma,\alpha)x^{\alpha, \mathbf{i}}\ (\textrm{mod}\
\mathscr{A}_{\langle \hat{k}-1 \rangle}).
\end{equation}

For all $r \in \overline{1,n-1}$, repeat the procedure of the proof
of \emph{Claim 2}. It is clear that \emph{Claim 3} holds.

Hence, by the minimality of $\iota (u)$, \emph{Claim 3} and Lemma
\ref{lemma:1}, we get
\begin{equation}\label{equ:2.55}
\phi^{2}(\gamma,\alpha)=\phi^{2}(\gamma,\beta)\quad \textrm{whenever
} a_{\alpha, \mathbf{i}}a_{\beta, \mathbf{j}} \not=0,
\end{equation}
namely,
\begin{equation}\label{equ:2.56}
\phi(\gamma,\alpha)=\phi(\gamma,\beta) \textrm{ or }
\phi(\gamma,\alpha)=-\phi(\gamma,\beta)\quad \textrm{whenever }
a_{\alpha, \mathbf{i}}a_{\beta, \mathbf{j}} \not=0.
\end{equation}
Assume that there exist $\alpha, \beta \in \Gamma$ such that
\begin{equation}
\phi(\gamma,\alpha)\not=\phi(\gamma,\beta),\
\phi(\gamma,\alpha)=-\phi(\gamma,\beta) \textrm{ and } a_{\alpha,
\mathbf{i}}a_{\beta, \mathbf{j}} \not= 0.
\end{equation}
Then we obtain $\phi(\gamma,\alpha)=-\phi(\gamma,\beta)\not= 0$.

For any $a_{\alpha, \mathbf{i}} \not = 0 $, if there exists
$q\in\overline{3,2n}$ such that $\varphi_{q}(\alpha)\not=0$ and
$\varphi_{q'}\not\equiv 0$. Picking $ \tau' \in
\textrm{Rad}_{\phi}\cap(\cap _{q' \not= p \in \overline{1,2n}}
\textrm{ker}_{\varphi_p} )\backslash \textrm{ker}_{\varphi_{q'}}$,
we have
\begin{eqnarray}
&&\mathcal{I}_{[0]}\cap \mathscr{A}_{\langle \hat{k} \rangle} \owns
[x^{\alpha+\tau', \mathbf{0}}, u]\nonumber\\
 & \equiv & \sum_{\substack{(-\varepsilon, \mathbf{0})\not =(\theta,
\mathbf{l})\in
\textrm{ker}_{\varphi_{2}} \times \boldsymbol{\mathcal{J}}\\
l_{2}=0, |\mathbf{l}|= \hat{k}}} a_{\theta, \mathbf{l}}
\epsilon_{q'} \varphi_{q'}(\tau')\varphi_{q}(\theta)
x^{\theta+\alpha+\tau'+\sigma_{q}+\varepsilon, \mathbf{l}} \
(\textrm{mod}\ \mathscr{A}_{\langle \hat{k}-1
\rangle}).\label{equ:2.57}
\end{eqnarray}
Since $2\alpha+\tau'+\sigma_{q}+\varepsilon\not=-\varepsilon$ due to
$\phi(\gamma,2\alpha+\tau'+\sigma_{q}+\varepsilon)
=2\phi(\gamma,\alpha) \not=0$, and $ a_{\alpha, \mathbf{i}}
\epsilon_{q'} \varphi_{q'}(\tau')\varphi_{q}(\alpha)\not=0, $ we
have $[x^{\alpha+\tau', \mathbf{0}}, u] \in (\mathcal{I}_{[0]}\cap
\mathscr{A}_{\langle \hat{k} \rangle}) \backslash
\mathbb{F}x^{-\varepsilon,\mathbf{0}}$. Then the minimality of
$\iota(u)$ implies that $\iota([x^{\alpha+\tau', \mathbf{0}},
u])=\iota(u)$. But we have
\begin{equation}
\phi(\gamma,2\alpha+\tau'+\sigma_{q}+\varepsilon)\not=\phi(\gamma,\beta+\alpha+\tau'+\sigma_{q}+\varepsilon)
\end{equation}
and
\begin{equation}
\phi(\gamma,2\alpha+\tau'+\sigma_{q}+\varepsilon)+\phi(\gamma,\beta+\alpha+\tau'+\sigma_{q}+\varepsilon)=2\phi(\gamma,\alpha)\not=0,
\end{equation}
which contradicts (\ref{equ:2.56}) if we replace $u$ by
$[x^{\alpha+\tau', \mathbf{0}}, u]$. Thus the first equation in
(\ref{equ:2.56}) holds, and so does (\ref{equ:2.35}).

For any $a_{\alpha, \mathbf{i}} \not = 0 $, if either
$\varphi_{p}(\alpha)=0$ or $\varphi_{p'}\equiv 0$ for
$p\in\overline{3,2n}$, and there does exist $q\in\overline{3,2n}$
such that $\varphi_{q}(\alpha)\not=0$. Fixing $p$, if
$\varphi_{p}(\alpha)\not=0$, we have $\varphi_{p'}\equiv 0$. Then
$\mathcal {J}_{p'} = \mathbb{N}$ by (\ref{equ:1.1}). By
(\ref{equ:2.32}), $a_{\beta, \mathbf{j}} \not = 0$ implies that
$j_{p'} = 0$. In the following, we would like to rewrite $u'$ in
(\ref{equ:2.12}) as
\begin{equation}\label{equ:2.72}
u'=  \sum_{\substack{(-\varepsilon, \mathbf{0})\not=(\beta,
\mathbf{l})
\in\textrm{ker}_{\varphi_{2}}\times\boldsymbol{\mathcal{J}}\\
l_{2}=0, |\mathbf{l}|= \hat{k}-1}} c_{\beta, \mathbf{l}} x^{\beta,
\mathbf{l}}+ u'' ,
\end{equation}
where $\ c_{\beta, \mathbf{l}}\in\mathbb{F}$ and $u'' \in
\mathscr{A}_{\langle \hat{k}-2 \rangle}+\mathbb{F}x^{-\varepsilon,
\mathbf{0}}$. Thus we get
\begin{eqnarray}
&&\mathcal{I}_{[0]}\cap \mathscr{A}_{\langle \hat{k}-1 \rangle}
\owns
[x^{\alpha, \mathbf{0}}, u]\nonumber\\
 & \equiv & \sum_{\substack{(-\varepsilon, \mathbf{0})\not =(\theta',
\mathbf{l'})\in
\textrm{ker}_{\varphi_{2}} \times \boldsymbol{\mathcal{J}}\\
l'_{2}=0, |\mathbf{l'}|= \hat{k}-1}} c_{\theta', \mathbf{l'}}
\phi(\alpha,\theta') x^{\theta'+\alpha+\varepsilon, \mathbf{l'}}
\;\; (\textrm{mod}\ \mathscr{A}_{\langle \hat{k}-2 \rangle}).
\end{eqnarray}
For any $c_{\theta', \mathbf{l'}}\not=0$, if
$\phi(\alpha,\theta')\not=0$, we have
$\theta'+\alpha+\varepsilon\not=-\varepsilon$ due to
$\phi(\alpha,\theta'+\alpha+\varepsilon)\not=0$. Then $[x^{\alpha,
\mathbf{0}}, u]\in (\mathcal{I}_{[0]}\cap \mathscr{A}_{\langle
\hat{k}-1 \rangle})\backslash
\mathbb{F}x^{-\varepsilon,\mathbf{0}}$, which contradicts the
minimality of $\hat{k}$ in (\ref{equ:2.11}). Therefore
$\phi(\alpha,\theta')=0$ whenever $c_{\theta', \mathbf{l'}}\not=0$.
Hence we have
\begin{eqnarray}
&&\mathcal{I}_{[0]}\cap \mathscr{A}_{\langle \hat{k} \rangle} \owns
[x^{\alpha, 1_{[q']}}, u]\nonumber\\
 & \equiv & \sum_{\substack{(-\varepsilon, \mathbf{0})\not =(\theta,
\mathbf{l})\in
\textrm{ker}_{\varphi_{2}} \times \boldsymbol{\mathcal{J}}\\
l_{2}=0, |\mathbf{l}|= \hat{k}}} a_{\theta, \mathbf{l}}
\epsilon_{q'} \varphi_{q}(\theta)
x^{\theta+\alpha+\sigma_{q}+\varepsilon, \mathbf{l}}
\;\;(\textrm{mod}\ \mathscr{A}_{\langle \hat{k}-1
\rangle}).\label{equ:2.58}
\end{eqnarray}
Since $2\alpha+\sigma_{q}+\varepsilon\not=-\varepsilon$ due to
$\phi(\gamma,2\alpha+\sigma_{q}+\varepsilon) \not=0$, and $
a_{\alpha, \mathbf{i}} \epsilon_{q'} \varphi_{q}(\alpha)\not=0,$ we
have $[x^{\alpha, 1_{[q']}}, u] \in (\mathcal{I}_{[0]}\cap
\mathscr{A}_{\langle \hat{k} \rangle}) \backslash
\mathbb{F}x^{-\varepsilon,\mathbf{0}}$. Then the minimality of
$\iota(u)$ implies that $\iota([x^{\alpha, 1_{[q']}}, u])=\iota(u)$.
But we have
\begin{equation}
\phi(\gamma,2\alpha+\sigma_{q}+\varepsilon)\not=\phi(\gamma,\beta+\alpha+\sigma_{q}+\varepsilon)
\end{equation}
and
\begin{equation}
\phi(\gamma,2\alpha+\sigma_{q}+\varepsilon)+\phi(\gamma,\beta+\alpha+\sigma_{q}+\varepsilon)=2\phi(\gamma,\alpha)\not=0,
\end{equation}
which contradicts (\ref{equ:2.56}) if we replace $u$ by $[x^{\alpha,
1_{[q']}}, u]$. Therefore the first equation in (\ref{equ:2.56})
holds, and so does (\ref{equ:2.35}).

For any $a_{\alpha, \mathbf{i}} \not = 0 $, if
$\varphi_{q}(\alpha)=0$ for all $q\in\overline{3,2n}$, we therefore
have $\varphi_{1}(\alpha)+1 \not=0$ by (\ref{equ:2.41}). Moreover,
if $\varphi_{2} \not\equiv 0$, picking $ \tau \in
\textrm{Rad}_{\phi}\cap(\cap _{2 \not= p \in \overline{1,2n}}
\textrm{ker}_{\varphi_p} )\backslash \textrm{ker}_{\varphi_{2}}$, we
have
\begin{eqnarray}
&&\mathcal{I}_{[0]}\cap \mathscr{A}_{\langle \hat{k} \rangle} \owns
[x^{-\tau, \mathbf{0}},[x^{\alpha+\tau, \mathbf{0}}, u]]\nonumber\\
 & \equiv & \sum_{\substack{(-\varepsilon, \mathbf{0})\not =(\theta,
\mathbf{l})\in
\textrm{ker}_{\varphi_{2}} \times \boldsymbol{\mathcal{J}}\\
l_{2}=0, |\mathbf{l}|= \hat{k}}} -2 a_{\theta, \mathbf{l}}
\varphi_{2}^{2}(\tau)(\varphi_{1}(\alpha)+1)^{2} x^{\theta+\alpha,
\mathbf{l}} \ (\textrm{mod}\ \mathscr{A}_{\langle \hat{k}-1
\rangle}).\label{equ:2.59}
\end{eqnarray}
Since $2\alpha \not=-\varepsilon$ due to $\phi(\gamma,2\alpha)
=2\phi(\gamma,\alpha) \not=0$, and $ a_{\alpha, \mathbf{i}}
\varphi_{2}^{2}(\tau)(\varphi_{1}(\alpha)+1)^{2} \not=0, $ we have
$[x^{-\tau, \mathbf{0}},[x^{\alpha+\tau, \mathbf{0}}, u]] \in
(\mathcal{I}_{[0]}\cap \mathscr{A}_{\langle \hat{k} \rangle})
\backslash \mathbb{F}x^{-\varepsilon,\mathbf{0}}$. Then the
minimality of $\iota(u)$ implies that $\iota([x^{-\tau,
\mathbf{0}},[x^{\alpha+\tau, \mathbf{0}}, u]])=\iota(u)$. But we
have
\begin{equation}
\phi(\gamma,2\alpha)\not=\phi(\gamma,\beta+\alpha) \textrm{ and }
\phi(\gamma,2\alpha)+\phi(\gamma,\beta+\alpha)=2\phi(\gamma,\alpha)\not=0,
\end{equation}
which contradicts (\ref{equ:2.56}) if we replace $u$ by $[x^{-\tau,
\mathbf{0}},[x^{\alpha+\tau, \mathbf{0}}, u]]$. Thus the first
equation in (\ref{equ:2.56}) holds, and so does (\ref{equ:2.35}).
Suppose that $\varphi_{2} \equiv 0$, then $\mathcal {J}_{2} =
\mathbb{N}$ by (\ref{equ:1.1}). By (\ref{equ:4.2}) and
(\ref{equ:2.72}) we have
\begin{eqnarray}
&&\mathcal{I}_{[0]}\cap \mathscr{A}_{\langle \hat{k}-1 \rangle}
\owns
[x^{\alpha, \mathbf{0}}, u]\nonumber\\
 & \equiv & \sum_{\substack{(-\varepsilon, \mathbf{0})\not =(\theta',
\mathbf{l'})\in
\textrm{ker}_{\varphi_{2}} \times \boldsymbol{\mathcal{J}}\\
l'_{2}=0, |\mathbf{l'}|= \hat{k}-1}} c_{\theta', \mathbf{l'}}
\phi(\alpha,\theta') x^{\theta'+\alpha+\varepsilon, \mathbf{l'}} \
(\textrm{mod}\ \mathscr{A}_{\langle \hat{k}-2 \rangle}).
\end{eqnarray}
If $\phi(\alpha,\theta')\not=0$ for some $c_{\theta',
\mathbf{l'}}\not=0$, we have
$\theta'+\alpha+\varepsilon\not=-\varepsilon$ due to
$\phi(\alpha,\theta'+\alpha+\varepsilon)\not=0$. Then we obtain that
$[x^{\alpha, \mathbf{0}}, u]\in (\mathcal{I}_{[0]}\cap
\mathscr{A}_{\langle \hat{k}-1 \rangle})\backslash
\mathbb{F}x^{-\varepsilon,\mathbf{0}}$, which contradicts the
minimality of $\hat{k}$ in (\ref{equ:2.11}). Therefore
$\phi(\alpha,\theta')=0$ whenever $c_{\theta', \mathbf{l'}}\not=0$.
Similarly, we have $\phi(\alpha,\theta)=0$ for any nonzero terms
$c_{\theta, \mathbf{l}}x^{\theta, \mathbf{l}}$ of $u$ with
$|\mathbf{l}|\leq \hat{k}-2$. Thus we get
\begin{eqnarray}
&&\mathcal{I}_{[0]}\cap \mathscr{A}_{\langle \hat{k} \rangle} \owns
[x^{\alpha, 1_{[2]}}, u]\nonumber\\
 & \equiv & \sum_{\substack{(-\varepsilon, \mathbf{0})\not =(\theta,
\mathbf{l})\in
\textrm{ker}_{\varphi_{2}} \times \boldsymbol{\mathcal{J}}\\
l_{2}=0, |\mathbf{l}|= \hat{k}}} -a_{\theta, \mathbf{l}}
(\varphi_{1}(\theta)+1) x^{\theta+\alpha, \mathbf{l}} \
(\textrm{mod}\ \mathscr{A}_{\langle \hat{k}-1
\rangle}).\label{equ:2.60}
\end{eqnarray}
Since $2\alpha\not=-\varepsilon$ due to $\phi(\gamma,2\alpha)
=2\phi(\gamma,\alpha) \not=0$, and $ a_{\alpha, \mathbf{i}}
(\varphi_{1}(\alpha)+1)\not=0,$ we have $[x^{\alpha, 1_{[2]}}, u]
\in (\mathcal{I}_{[0]}\cap \mathscr{A}_{\langle \hat{k} \rangle})
\backslash \mathbb{F}x^{-\varepsilon,\mathbf{0}}$. Then the
minimality of $\iota(u)$ implies that $\iota([x^{\alpha, 1_{[2]}},
u])=\iota(u)$. But we have
\begin{equation}
\phi(\gamma,2\alpha)\not=\phi(\gamma,\beta+\alpha) \textrm{ and }
\phi(\gamma,2\alpha)+\phi(\gamma,\beta+\alpha)=2\phi(\gamma,\alpha)\not=0,
\end{equation}
which contradicts (\ref{equ:2.56}) if we replace $u$ by $[x^{\alpha,
1_{[2]}}, u]$. Thus the first equation in (\ref{equ:2.56}) holds,
and so does (\ref{equ:2.35}). Therefore we have proved that
\begin{equation}\label{equ:2.61}
\phi(\gamma, \alpha)=\phi(\gamma, \beta)\textrm{ for any }\gamma \in
\Gamma,\quad \textrm{ whenever } a_{\alpha, \mathbf{i}}a_{\beta,
\mathbf{j}}\not=0.
\end{equation}

Hence by (\ref{equ:1.3}), (\ref{equ:2.34}) and (\ref{equ:2.61}), we
get
\begin{equation}\label{equ:2.62}
\alpha = \beta \qquad\textrm{ whenever } a_{\alpha,
\mathbf{i}}a_{\beta, \mathbf{j}}\not=0.
\end{equation}

\noindent{\bf \emph{Step 2. $1\in\mathcal {I}.$}}\vspace{0.3cm}

Let $a_{\alpha, \mathbf{i}}\not=0$ be fixed. We will give the proof
in two cases.

\noindent{\emph{Case 1}}. $\hat{k}=0$.

In this case, $x^{\alpha, \mathbf{0}} \in \mathcal{I}_{[0]}$ for
some $\alpha \in \textrm{ker}_{\varphi_{2}} \backslash
\{-\varepsilon\}$. Assume that $\alpha \not=0$. Since $\alpha \not=
-\varepsilon$, we have $\varphi_{1} (\alpha)+1\not=0$ or
$\varphi_{q} (\alpha)\not=0$ for some $q\in \overline{3,2n}$, or
$\alpha \not\in \textrm{Rad}_{\phi}$.

\noindent{\emph{Subcase 1.1}}. $\varphi_{1} (\alpha)+1\not=0$.

If $\varphi_{2} \not\equiv 0$, picking $\tau \in \textrm{Rad}_{\phi}
\cap(\cap _{2 \not= p \in \overline{1,2n}} \textrm{ker}_{\varphi_p})
\backslash \textrm{ker}_{\varphi_2}$, we have
\begin{equation}\label{equ:2.63}
[[x^{\alpha, \mathbf{0}}, x^{-\alpha+\tau, \mathbf{0}}], x^{-\tau,
\mathbf{0}}]=-2 \varphi_{2}^{2}(\tau) (\varphi_{1}(\alpha)+1) x^{0,
\mathbf{0}} \in \mathcal{I}.
\end{equation}
Thus we get $1 \in \mathcal{I}$. Otherwise $\varphi_{2}\equiv 0$, we
have $\mathcal{J}_{2}=\mathbb{N}$. Then we get
\begin{equation}\label{equ:2.64}
[x^{\alpha, \mathbf{0}}, x^{-\alpha,
1_{[2]}}]=(\varphi_{1}(\alpha)+1) x^{0, \mathbf{0}} \in \mathcal{I}.
\end{equation}
So we obtain $1 \in \mathcal{I}$.

\noindent{\emph{Subcase 1.2}}. $\varphi_{1} (\alpha)+1=0$ and
$\varphi_{q} (\alpha)\not=0$ for some $q\in \overline{3,2n}$.

If $\varphi_{q'} \not\equiv 0$, picking $\tau' \in
\textrm{Rad}_{\phi} \cap (\cap _{q' \not= p \in \overline{1,2n}}
\textrm{ker}_{\varphi_p}) \backslash \textrm{ker}_{\varphi_{q'}}$,
we have
\begin{equation}\label{equ:2.65}
[x^{\alpha, \mathbf{0}}, x^{\tau', \mathbf{0}}]=\epsilon_{q}
\varphi_{q'}(\tau') \varphi_{q}(\alpha)
x^{\alpha+\tau'+\sigma_{q}+\varepsilon, \mathbf{0}} \in \mathcal{I}.
\end{equation}
Then $x^{\alpha+\tau'+\sigma_{q}+\varepsilon, \mathbf{0}}\ \in
\mathcal{I}$. Since $\varphi_{1}
(\alpha+\tau'+\sigma_{q}+\varepsilon)+1\not=0$ and $\varphi_{2}
(\alpha+\tau'+\sigma_{q}+\varepsilon)=0$, we go back to
\emph{Subcase 1.1} with $x^{\alpha, \mathbf{0}}$ replaced by
$x^{\alpha+\tau'+\sigma_{q}+\varepsilon, \mathbf{0}}$. Assume
$\varphi_{q'} \equiv 0$, then $\mathcal{J}_{q'}=\mathbb{N}$. So we
have
\begin{equation}\label{equ:2.66}
[x^{\alpha, \mathbf{0}}, x^{0, 1_{[q']}}]=\epsilon_{q}
\varphi_{q}(\alpha) x^{\alpha+\sigma_{q}+\varepsilon, \mathbf{0}}
\in \mathcal{I},
\end{equation}
which implies $x^{\alpha+\sigma_{q}+\varepsilon, \mathbf{0}}\ \in
\mathcal{I}$. Since $\varphi_{1}
(\alpha+\sigma_{q}+\varepsilon)+1\not=0$ and $\varphi_{2}
(\alpha+\sigma_{q}+\varepsilon)=0$, we go back to \emph{Subcase 1.1}
with $x^{\alpha, \mathbf{0}}$ replaced by
$x^{\alpha+\sigma_{q}+\varepsilon, \mathbf{0}}$.

\noindent{\emph{Subcase 1.3}}. $\varphi_{1} (\alpha)+1=0$ and
$\varphi_{q} (\alpha)=0$ for all $q\in \overline{3,2n}$.

Since $\alpha\not=-\varepsilon$, we have
$\alpha\not\in\textrm{Rad}_{\phi}$. Namely, there exists
$\beta\in\Gamma$ such that $\phi(\alpha,\beta)\not=0$. Fixing any
such $\beta$, we have
\begin{equation}\label{equ:2.67}
[x^{\alpha, \mathbf{0}}, x^{\beta-\alpha-\varepsilon,
\mathbf{0}}]=\phi(\alpha,\beta) x^{\beta, \mathbf{0}} \in
\mathcal{I},
\end{equation}
which implies $x^{\beta, \mathbf{0}}\ \in \mathcal{I}$. If
$\varphi_{1} (\beta)+1\not=0$ and $\varphi_{2}(\beta)=0$, we go back
to \emph{Subcase 1.1} with $x^{\alpha, \mathbf{0}}$ replaced by
$x^{\beta, \mathbf{0}}$. If $\varphi_{1}(\beta)+1 \not=0$ and
$\varphi_{2}(\beta)\not=0$, we have
\begin{equation}\label{equ:2.68}
[x^{\beta, \mathbf{0}}, x^{-\beta, \mathbf{0}}]=-2
\varphi_{2}(\beta) x^{0, \mathbf{0}} \in \mathcal{I}.
\end{equation}
Then we get $1 \in \mathcal{I}$. If $\varphi_{1}(\beta)+1 =0$, we
have $\varphi_{1}(\beta+\varepsilon)+1\not=0$ and
$\phi(\alpha,\beta+\varepsilon)\not=0$. We go back to
(\ref{equ:2.67}) and (\ref{equ:2.68}) with $\beta$ replaced by
$\beta+\varepsilon$.

\noindent{\emph{Case 2}}. $\hat{k}>0$.

We have obtained that
\begin{equation}\label{equ:2.71}
u= \sum_{\substack{\mathbf{i}\in
\boldsymbol{\mathcal{J}}\\
i_{2}=0, |\mathbf{i}|= \hat{k}}} a_{\alpha, \mathbf{i}} x^{\alpha,
\mathbf{i}}+u' \ \in (\mathcal{I}_{[0]}\cap \mathscr{A}_{\langle
\hat{k} \rangle}) \backslash \mathbb{F}x^{-\varepsilon, \mathbf{0}},
\end{equation}
where $\alpha\in\textrm{ker}_{\varphi_{2}} \backslash
\{-\varepsilon\}$, and $u' \in \mathscr{A}_{\langle \hat{k}-1
\rangle}+\mathbb{F}x^{-\varepsilon, \mathbf{0}}$. In the following,
we will use $u'$ frequently in its specific form as in (\ref{equ:2.72}).

\noindent{\emph{Subcase 2.1}}. There exists $q\in\overline{3,2n}$
and $a_{\alpha, \mathbf{i}}\not=0$ such that $i_{q}>0$ and
$\varphi_{q} (\alpha)=0$.

If $\varphi_{q'}\not\equiv 0$, choose $\tau' \in
\textrm{Rad}_{\phi}\cap(\cap _{q' \not= p \in \overline{1,2n}}
\textrm{ker}_{\varphi_p} )\backslash \textrm{ker}_{\varphi_{q'}}$
such that
\begin{equation}\label{equ:2.102}
\varphi_{q'}(\alpha+\tau'+\sigma_{q})\not=0.
\end{equation}
Then we have
\begin{eqnarray}
[x^{\tau', \mathbf{0}}, u] & \equiv & (\sum_{\substack{\mathbf{i}\in
\boldsymbol{\mathcal{J}}\\
i_{2}=0, |\mathbf{i}|= \hat{k}}} a_{\alpha, \mathbf{i}}
\epsilon_{q'} i_{q} \varphi_{q'} (\tau') x^{\alpha+\tau'+\sigma_{q}+
\varepsilon, \mathbf{i}-1_{[q]}} \nonumber\\
& &+ \sum_{\substack{\alpha+\tau'+\sigma_{q}+ \varepsilon \not
=\beta
\in\Gamma, \mathbf{l} \in \boldsymbol{\mathcal{J}}\\
l_{2}=0, |\mathbf{l}|= \hat{k}-1}} c_{\beta, \mathbf{l}}' x^{\beta,
\mathbf{l}}) \
 (\textrm{mod}\ \mathscr{A}_{\langle \hat{k}-2 \rangle})\label{equ:2.73}
\end{eqnarray}
with $c_{\beta, \mathbf{l}}'\in\mathbb{F}$. Since
$\varphi_{q'}(\alpha+\tau'+\sigma_{q}+\varepsilon)\not=0$ by
(\ref{equ:2.102}), we have
$\alpha+\tau'+\sigma_{q}+\varepsilon\not=-\varepsilon$. Then
$[x^{\tau', \mathbf{0}}, u] \in (\mathcal{I}_{[0]}\cap
\mathscr{A}_{\langle \hat{k}-1 \rangle}) \backslash
\mathbb{F}x^{-\varepsilon, \mathbf{0}}$, which contradicts the
minimality of $\hat{k}$. Assume that $\varphi_{q'}\equiv 0$, then
$\mathcal{J}_{q'}=\mathbb{N}$. Picking $\tau \in
\textrm{Rad}_{\phi}\cap(\cap _{p\in \overline{2,2n}}
\textrm{ker}_{\varphi_p}) \backslash \textrm{ker}_{\varphi_{1}}$
such that $\varphi_{1}(\alpha+\tau+ \varepsilon)+1\not=0$, we have
\begin{eqnarray}
[x^{\tau, 1_{[q']}}, u] &\equiv & (\sum_{\substack{\mathbf{i}\in
\boldsymbol{\mathcal{J}}\\
i_{2}=0, |\mathbf{i}|= \hat{k}}} a_{\alpha, \mathbf{i}}
\epsilon_{q'} i_{q} x^{\alpha+\tau+\sigma_{q}+ \varepsilon,
\mathbf{i}-1_{[q]}}\nonumber \\
& & + \sum_{\substack{\alpha+\tau+\sigma_{q}+ \varepsilon \not
=\beta
\in\Gamma, \mathbf{l} \in \boldsymbol{\mathcal{J}}\\
l_{2}=0, |\mathbf{l}|= \hat{k}-1}} c_{\beta, \mathbf{l}}' x^{\beta,
\mathbf{l}} )\
 (\textrm{mod}\ \mathscr{A}_{\langle \hat{k}-2
\rangle})\label{equ:2.74}
\end{eqnarray}
with $c_{\beta, \mathbf{l}}'\in\mathbb{F}$. Since
$\varphi_{1}(\alpha+\tau+\sigma_{q}+ \varepsilon)+1\not=0$, we have
$\alpha+\tau+\sigma_{q}+ \varepsilon\not=-\varepsilon$. Then
$[x^{\tau, 1{[q']}}, u] \in (\mathcal{I}_{[0]}\cap
\mathscr{A}_{\langle \hat{k}-1 \rangle}) \backslash
\mathbb{F}x^{-\varepsilon, \mathbf{0}}$, which contradicts the
minimality of $\hat{k}$.

\noindent{\emph{Subcase 2.2}}. There exists $q\in\overline{3,2n}$
and $a_{\alpha, \mathbf{i}}\not=0$ such that $i_{q}>0$ and
$\varphi_{q} (\alpha)\not=0$.

If $\varphi_{q'}\not\equiv 0$, choose $\tau' \in
\textrm{Rad}_{\phi}\cap(\cap _{q' \not= p\in \overline{1,2n}}
\textrm{ker}_{\varphi_p}) \backslash \textrm{ker}_{\varphi_{q'}}$
such that
\begin{equation}
\varphi_{q} (\alpha)\varphi_{q'} (\tau'-\sigma_{q})+\varphi_{q'}
(\alpha)\varphi_{q} (\sigma_{q})\not=0.
\end{equation}
Then we have
\begin{eqnarray}
 (\mathcal{I}_{[0]}\cap \mathscr{A}_{\langle \hat{k}
\rangle}) \backslash \mathbb{F}x^{-\varepsilon, \mathbf{0}} & \owns
& [x^{\tau'-\alpha-\sigma_{q}-\varepsilon, \mathbf{0}}, u]  \equiv
\sum_{\substack{\mathbf{i}\in
\boldsymbol{\mathcal{J}}\\
i_{2}=0, |\mathbf{i}|= \hat{k}}} a_{\alpha, \mathbf{i}}
\epsilon_{q'} \big( \varphi_{q} (\alpha)\varphi_{q'}
(\tau'-\sigma_{q})
\nonumber\\
& & +\varphi_{q'} (\alpha)\varphi_{q} (\sigma_{q})\big) x^{\tau',
\mathbf{i}}\;\; (\textrm{mod}\ \mathscr{A}_{\langle \hat{k}-1
\rangle}).\label{equ:2.75}
\end{eqnarray}
Since $\varphi_{q}(\tau')=0$, we go back to \emph{Subcase 2.1} with
$u$ replaced by $[x^{\tau'-\alpha-\sigma_{q}-\varepsilon,
\mathbf{0}}, u]$. Assume that $\varphi_{q'}\equiv 0$, by
(\ref{equ:1.1}) we have $\mathcal{J}_{q'}=\mathbb{N}$. Then we get
\begin{eqnarray}\label{equ:2.76}
&& \mathcal{I}_{[0]}\cap \mathscr{A}_{\langle \hat{k}-1
\rangle} \owns [x^{-\alpha-\sigma_{q}, \mathbf{0}}, u] \nonumber\\
& \equiv & \sum_{\substack{\mathbf{i}\in
\boldsymbol{\mathcal{J}}\\
i_{2}=0, |\mathbf{i}|= \hat{k}}} \sum_{s=1}^{n-1} a_{\alpha,
\mathbf{i}}
\big(i_{2s+1}\varphi_{2s+2}(\alpha+\sigma_{q})x^{\sigma_{2s+1}+\varepsilon-\sigma_{q},\mathbf{i}-1_{[2s+1]}}
- i_{2s+2}\varphi_{2s+1}(\alpha+\sigma_{q})\cdot \nonumber\\
& &
x^{\sigma_{2s+1}+\varepsilon-\sigma_{q},\mathbf{i}-1_{[2s+2]}}\big)
+\sum_{\substack{(-\varepsilon, \mathbf{0})\not=(\beta, \mathbf{l})
\in\textrm{ker}_{\varphi_{2}}\times\boldsymbol{\mathcal{J}}\\
l_{2}=0, |\mathbf{l}|= \hat{k}-1}} c_{\beta,
\mathbf{l}} \big(\sum_{s=1}^{n-1}(\varphi_{2s+1}(\beta)\varphi_{2s+2}(\alpha+\sigma_{q})\nonumber\\
& &-\varphi_{2s+1}(\alpha+\sigma_{q})\varphi_{2s+2}(\beta))
x^{\beta+\sigma_{2s+1}+\varepsilon-\alpha-\sigma_{q},\mathbf{l}}\nonumber\\
 & & -\phi(\alpha,\beta)x^{\beta+\varepsilon-\alpha-\sigma_{q},\mathbf{l}} \big)\;\;(\textrm{mod}\
\mathscr{A}_{\langle \hat{k}-2 \rangle})
\end{eqnarray}
and
\begin{eqnarray}\label{equ:2.77}
&& \mathcal{I}_{[0]}\cap \mathscr{A}_{\langle \hat{k}-1
\rangle} \owns [x^{-\alpha-\sigma_{q}-\varepsilon, \mathbf{0}}, u] \nonumber\\
& \equiv & \sum_{\substack{\mathbf{i}\in
\boldsymbol{\mathcal{J}}\\
i_{2}=0, |\mathbf{i}|= \hat{k}}} \sum_{s=1}^{n-1} a_{\alpha,
\mathbf{i}}
(i_{2s+1}\varphi_{2s+2}(\alpha+\sigma_{q})x^{\sigma_{2s+1}-\sigma_{q},\mathbf{i}-1_{[2s+1]}}
-
i_{2s+2}\varphi_{2s+1}(\alpha+\sigma_{q})\nonumber\\
& & \cdot x^{\sigma_{2s+1}-\sigma_{q},\mathbf{i}-1_{[2s+2]}})
+\sum_{\substack{(-\varepsilon, \mathbf{0})\not=(\beta, \mathbf{l})
\in\textrm{ker}_{\varphi_{2}}\times\boldsymbol{\mathcal{J}}\\
l_{2}=0, |\mathbf{l}|= \hat{k}-1}} c_{\beta,
\mathbf{l}}(\sum_{s=1}^{n-1}(\varphi_{2s+1}(\beta)\varphi_{2s+2}(\alpha+\sigma_{q})\nonumber\\
& &-\varphi_{2s+1}(\alpha+\sigma_{q})\varphi_{2s+2}(\beta))
x^{\beta+\sigma_{2s+1}-\alpha-\sigma_{q},\mathbf{l}}\nonumber\\
 & & -\phi(\alpha,\beta)x^{\beta-\alpha-\sigma_{q},\mathbf{l}})\;\;(\textrm{mod}\
\mathscr{A}_{\langle \hat{k}-2 \rangle}).
\end{eqnarray}
Since $\sigma_{2s+1}-\sigma_{q}\not=-\varepsilon$ for all
$s\in\overline{1,n-1}$, avoiding contradiction to the minimality of
$\hat{k}$, the coefficients corresponding to the terms
$x^{\sigma_{2s+1}-\sigma_{q},\mathbf{i}-1_{[2s+1]}}$ and
$x^{\sigma_{2s+1}-\sigma_{q},\mathbf{i}-1_{[2s+2]}}$ for all
$s\in\overline{1,n-1}$ in (\ref{equ:2.77}) are zeros. For any
$s\in\overline{1,n-1}$ , we have
\begin{equation}
\beta+\sigma_{2s+1}+\varepsilon-\alpha-\sigma_{q}\not=-\varepsilon
\textrm{ or }
\beta+\sigma_{2s+1}-\alpha-\sigma_{q}\not=-\varepsilon.
\end{equation}
The coefficient corresponding to the term
$x^{\beta+\sigma_{2s+1}+\varepsilon-\alpha-\sigma_{q},\mathbf{l}}$
in (\ref{equ:2.76}) and the coefficient corresponding to the term
$x^{\beta+\sigma_{2s+1}-\alpha-\sigma_{q},\mathbf{l}}$ in
(\ref{equ:2.77}) are the same. To avoid contradiction to the
minimality of $\hat{k}$, we get that these coefficients equal zero.
Similarly, we have
\begin{equation}
\beta+\varepsilon-\alpha-\sigma_{q}\not=-\varepsilon \textrm{ or }
\beta-\alpha-\sigma_{q}\not=-\varepsilon.
\end{equation}
The coefficient corresponding to the term
$x^{\beta+\varepsilon-\alpha-\sigma_{q},\mathbf{l}}$ in
(\ref{equ:2.76}) and the coefficient corresponding to the term
$x^{\beta-\alpha-\sigma_{q},\mathbf{l}}$ in (\ref{equ:2.77}) are the
same. To avoid contradiction to the minimality of $\hat{k}$, we get
that these coefficients equal zero too. Using these results, we have
\begin{eqnarray}\label{equ:2.79}
&& \mathcal{I}_{[0]}\cap \mathscr{A}_{\langle \hat{k}
\rangle}\owns [x^{-\alpha-\sigma_{q}-\varepsilon, 1_{[q']}}, u] \nonumber\\
& \equiv & \sum_{\substack{\mathbf{i}\in
\boldsymbol{\mathcal{J}}\\
i_{2}=0, |\mathbf{i}|= \hat{k}}} a_{\alpha, \mathbf{i}}
\epsilon_{q'} \varphi_{q}(\alpha) x^{0, \mathbf{i}}\ (\textrm{mod}\
\mathscr{A}_{\langle \hat{k}-1 \rangle}).
\end{eqnarray}
Thus we go back to \emph{Subcase 2.1} with $u$ replaced by
$[x^{-\alpha-\sigma_{q}-\varepsilon, 1_{[q']}}, u]$.

\noindent{\emph{Subcase 2.3}}. For any $a_{\alpha,
\mathbf{i}}\not=0$, $i_{q}=0$ for all $q\in\overline{3,2n}$.

Since $i_{q}=0$ for all $q\in\overline{2,2n}$, we obtain
$i_{1}=\hat{k}$. Then we have
\begin{equation}\label{equ:2.80}
u= a_{\alpha, \hat{k}_{[1]}} x^{\alpha, \hat{k}_{[1]}} + u' \ \in
(\mathcal{I}_{[0]}\cap \mathscr{A}_{\langle \hat{k} \rangle})
\backslash \mathbb{F}x^{-\varepsilon, \mathbf{0}}.
\end{equation}

If $\varphi_{2} \not\equiv 0$, pick $\tau \in
\textrm{Rad}_{\phi}\cap (\cap_{2\not= p \in
\overline{1,2n}}\textrm{ker}_{\varphi_p})\backslash
\textrm{ker}_{\varphi_2}$. Moreover, if $\varphi_{1}(\alpha) +
1\not= 0$, we have
\begin{eqnarray}\label{equ:2.82}
&& (\mathcal{I}_{[0]}\cap \mathscr{A}_{\langle \hat{k} \rangle})
\backslash \mathbb{F}x^{-\varepsilon, \mathbf{0}} \owns
[x^{-\tau-\varepsilon, \mathbf{0}},[x^{\tau-\alpha,
\mathbf{0}},u]]\nonumber\\
& \equiv & -a_{\alpha, \hat{k}_{[1]}} \varphi_{2}^{2} (\tau)
(\varphi_{1}(\alpha ) +1) x^{-\varepsilon, \hat{k}_{[1]}}\
(\textrm{mod}\ \mathscr{A}_{\langle \hat{k}-1 \rangle}).
\end{eqnarray}
Since $\varphi_{1}(-\varepsilon) +1=0$, replacing $u$ by
$[x^{-\tau-\varepsilon, \mathbf{0}},[x^{\tau-\alpha,
\mathbf{0}},u]]$ and $\alpha$ by $-\varepsilon$, we can always
assume that $\varphi_{1}(\alpha) + 1 = 0$. Choose $\tau' \in
\textrm{Rad}_{\phi}\cap (\cap _{ p \in \overline{2,2n}}
\textrm{ker}_{\varphi_p}) \backslash \textrm{ker}_{\varphi_1}$ such
that $ \varphi_{1}(\tau') +1 \not=0$. Then we have
\begin{eqnarray}\label{equ:2.81}
[x^{-\tau+\tau', \mathbf{0}},[ x^{\tau, \mathbf{0}},u] ]& \equiv &
-a_{\alpha, \hat{k}_{[1]}} \hat{k} \varphi_{2}^{2} (\tau)
(\varphi_{1}(\tau') +1) x^{\alpha+\tau', (\hat{k}-1)_{[1]}} \nonumber \\
& & + \sum_{\substack{\alpha+\tau' \not =\beta
\in\Gamma, \mathbf{l} \in \boldsymbol{\mathcal{J}}\\
l_{2}=0, |\mathbf{l}|= \hat{k}-1}} c_{\beta, \mathbf{l}}' x^{\beta,
\mathbf{l}} \ (\textrm{mod}\ \mathscr{A}_{\langle \hat{k}-2
\rangle})
\end{eqnarray}
with $c_{\beta, \mathbf{l}}'\in\mathbb{F}$. Since
$\varphi_{1}(\alpha+\tau')+1 \not=0$, we get
$\alpha+\tau'\not=-\varepsilon$. Then we have $[x^{-\tau+\tau',
\mathbf{0}},[ x^{\tau, \mathbf{0}},u] ] \in (\mathcal{I}_{[0]}\cap
\mathscr{A}_{\langle \hat{k}-1 \rangle}) \backslash
\mathbb{F}x^{-\varepsilon, \mathbf{0}}$, which contradicts the
minimality of $\hat{k}$.

Assume that $\varphi_{2} \equiv 0$, by (\ref{equ:1.1}) we have
$\mathcal{J}_{2}=\mathbb{N}$. Then we get
\begin{eqnarray}
&& \mathcal{I}_{[0]}\cap \mathscr{A}_{\langle \hat{k}-1
\rangle} \owns [x^{-\alpha, \mathbf{0}}, u] \nonumber\\
& \equiv & \sum_{\substack{(-\varepsilon, \mathbf{0})\not=(\beta,
\mathbf{l})
\in\textrm{ker}_{\varphi_{2}}\times\boldsymbol{\mathcal{J}}\\
l_{2}=0, |\mathbf{l}|= \hat{k}-1}} c_{\beta, \mathbf{l}} \big(
\sum_{s=1}^{n-1}(\varphi_{2s+1}(\beta)\varphi_{2s+2}(\alpha)
-\varphi_{2s+1}(\alpha)\varphi_{2s+2}(\beta))\cdot\nonumber\\
& &  x^{\beta+\sigma_{2s+1}+\varepsilon-\alpha,\mathbf{l}}
  -\phi(\alpha,\beta)x^{\beta+\varepsilon-\alpha,\mathbf{l}}\big)\;\;(\textrm{mod}\
\mathscr{A}_{\langle \hat{k}-2 \rangle} )\label{equ:2.84}
\end{eqnarray}
and
\begin{eqnarray}
&& \mathcal{I}_{[0]}\cap \mathscr{A}_{\langle \hat{k}-1
\rangle} \owns [x^{-\alpha-\varepsilon, \mathbf{0}}, u] \nonumber\\
& \equiv & \sum_{\substack{(-\varepsilon, \mathbf{0})\not=(\beta,
\mathbf{l})
\in\textrm{ker}_{\varphi_{2}}\times\boldsymbol{\mathcal{J}}\\
l_{2}=0, |\mathbf{l}|= \hat{k}-1}} c_{\beta, \mathbf{l}} \big(
\sum_{s=1}^{n-1}(\varphi_{2s+1}(\beta)\varphi_{2s+2}(\alpha)
-\varphi_{2s+1}(\alpha)\varphi_{2s+2}(\beta)) \cdot \nonumber\\
 & & x^{\beta+\sigma_{2s+1}-\alpha,\mathbf{l}}-\phi(\alpha,\beta)x^{\beta-\alpha,\mathbf{l}} \big)\;\;(\textrm{mod}\
\mathscr{A}_{\langle \hat{k}-2 \rangle}).\label{equ:2.85}
\end{eqnarray}
For any $s\in\overline{1,n-1}$, we have
\begin{equation}
\beta+\sigma_{2s+1}+\varepsilon-\alpha\not=-\varepsilon \textrm{ or
} \beta+\sigma_{2s+1}-\alpha\not=-\varepsilon.
\end{equation}
The coefficient corresponding to the term
$x^{\beta+\sigma_{2s+1}+\varepsilon-\alpha,\mathbf{l}}$ in
(\ref{equ:2.84}) and the coefficient corresponding to the term
$x^{\beta+\sigma_{2s+1}-\alpha,\mathbf{l}}$ in (\ref{equ:2.85}) are
the same. To avoid contradiction to the minimality of $\hat{k}$, we
get that these coefficients equal zero. Similarly, we have
\begin{equation}
\beta+\varepsilon-\alpha\not=-\varepsilon \textrm{ or }
\beta-\alpha\not=-\varepsilon.
\end{equation}
The coefficient corresponding to the term
$x^{\beta+\varepsilon-\alpha,\mathbf{l}}$ in (\ref{equ:2.84}) and
the coefficient corresponding to the term
$x^{\beta-\alpha,\mathbf{l}}$ in (\ref{equ:2.85}) are the same. To
avoid contradiction to the minimality of $\hat{k}$, we get that
these coefficients equal zero too. Using these results, we have
\begin{eqnarray}
&& \mathcal{I}_{[0]}\cap \mathscr{A}_{\langle \hat{k}
\rangle}\owns [x^{-\alpha-\varepsilon, 1_{[2]}}, u] \nonumber\\
& \equiv & - a_{\alpha, \hat{k}_{[1]}} (\varphi_{1}(\alpha)+1)
x^{-\varepsilon, \hat{k}_{[1]}}\ (\textrm{mod}\ \mathscr{A}_{\langle
\hat{k}-1 \rangle}).\label{equ:2.86}
\end{eqnarray}
If  $\varphi_{1}(\alpha) + 1 \not= 0$, by replacing $u$ by
$[x^{-\alpha-\varepsilon, 1_{[2]}}, u]$ and $\alpha$ by
$-\varepsilon$, we can always assume that $\varphi_{1}(\alpha) + 1 =
0$. Picking $\tau' \in \textrm{Rad}_{\phi}\cap (\cap _{ p \in
\overline{2,2n}} \textrm{ker}_{\varphi_p} )\backslash
\textrm{ker}_{\varphi_1}$, we have
\begin{equation}\label{equ:2.83}
[x^{\tau',1_{[2]}},u] \equiv -\hat{k} a_{\alpha, \hat{k}_{[1]}}
x^{\alpha+\tau', (\hat{k}-1)_{[1]}} +
\sum_{\substack{\alpha+\tau'\not =\beta
\in\Gamma, \mathbf{l} \in \boldsymbol{\mathcal{J}}\\
l_{2}=0, |\mathbf{l}|= \hat{k}-1}} c_{\beta, \mathbf{l}}' x^{\beta,
\mathbf{l}} \ (\textrm{mod}\ \mathscr{A}_{\langle \hat{k}-2
\rangle})
\end{equation}
with $c_{\beta, \mathbf{l}}'\in\mathbb{F}$. Since
$\varphi_{1}(\alpha+\tau')+1 = \varphi_{1}(\tau')\not=0$, we have
$\alpha+\tau'\not=-\varepsilon$. Then $[x^{\tau',1_{[2]}},u] \in
(\mathcal{I}_{[0]}\cap \mathscr{A}_{\langle \hat{k}-1 \rangle})
\backslash \mathbb{F}x^{-\varepsilon, \mathbf{0}}$, which
contradicts the minimality of $\hat{k}$.

So we always have $\hat{k}=0$ and $1 \in \mathcal{I}$.\vspace{0.3cm}

\noindent{\bf \emph{Step 3. The conclusion of Step 2 implies
$\mathcal{I}=\mathscr{A}$ if
$\boldsymbol{\mathcal{J}}\not=\{0\}$.}}\vspace{0.3cm}

Since
\begin{equation}\label{equ:2.87}
[1, x^{\beta, \mathbf{j}}]=\varphi_{2}(\beta)x^{\beta,
\mathbf{j}}+j_{2}x^{\beta, \mathbf{j}-1_{[2]}} \qquad\textrm{for}\
(\beta,\mathbf{j})\in\Gamma\times\boldsymbol{\mathcal{J}}.
\end{equation}
If $\mathcal{J}_{2} = \mathbb{N}$, we can prove
$\mathcal{I}=\mathscr{A}$ by induction on $j_{2}$. Assume that
$\mathcal{J}_{2} = \{0\}$, then $\varphi_{2}\not\equiv 0$ by
(\ref{equ:1.1}). Thus (\ref{equ:2.87}) implies
\begin{equation}\label{equ:2.88}
x^{\beta, \mathbf{j}}\in \mathcal{I}\quad \textrm{for}\ (\beta,
\mathbf{j})\in \Gamma\times\boldsymbol{\mathcal{J}}, \
\varphi_{2}(\beta)\not=0.
\end{equation}
Picking $ \tau \in \textrm{Rad}_{\phi}\cap(\cap _{2 \not= p \in
\overline{1,2n}} \textrm{ker}_{\varphi_p}) \backslash
\textrm{ker}_{\varphi_2}$, we get
$x^{\tau,\mathbf{0}}\in\mathcal{I}$. Since we have
\begin{equation}\label{equ:2.89}
\mathcal{I}\owns [x^{\tau,\mathbf{0}}, x^{-\tau+\beta,\mathbf{j}}]=
-\varphi_{2}(\tau)(\varphi_{1}(\beta)+2)x^{\beta,\mathbf{j}}-
j_{1}\varphi_{2}(\tau)x^{\beta,\mathbf{j}-1_{[1]}}
\end{equation}
for any $(\beta, \mathbf{j})\in
\textrm{ker}_{\varphi_{2}}\times\boldsymbol{\mathcal{J}}$. If
$\mathcal{J}_{1}=\mathbb{N}$, we can prove $\mathcal{I}
=\mathscr{A}$ by induction on $j_{1}$. Assume that
$\mathcal{J}_{1}=\{0\}$. Then (\ref{equ:2.89}) shows
\begin{equation}\label{equ:2.90}
x^{\beta, \mathbf{j}}\in \mathcal{I}\quad \textrm{for}\ (\beta,
\mathbf{j})\in
\textrm{ker}_{\varphi_{2}}\times\boldsymbol{\mathcal{J}},
\varphi_{1}(\beta)+2\not=0.
\end{equation}
If there exists $q\in\overline{3,2n}$ such that
$\mathcal{J}_{q}\not=\{0\}$. Moreover, if $\varphi_{q'} \not\equiv
0$, picking $ \tau \in \textrm{Rad}_{\phi} \cap(\cap _{q' \not= p
\in \overline{1,2n}} \textrm{ker}_{\varphi_p} )\backslash
\textrm{ker}_{\varphi_{q'}}$, we have $x^{\tau, \mathbf{0}} \in
\mathcal{I}$ due to $\varphi_{1}(\tau)+2\not=0$. Since we get
\begin{equation}\label{equ:2.91}
\mathcal{I}\owns [x^{\tau,\mathbf{0}},
x^{\beta-\tau-\sigma_{q}-\varepsilon,\mathbf{j}}]=\epsilon_{q'}(
\varphi_{q'}(\tau)\varphi_{q}(\beta-\sigma_{q})x^{\beta,\mathbf{j}}+
j_{q}\varphi_{q'}(\tau)x^{\beta,\mathbf{j}-1_{[q]}})
\end{equation}
for $\beta\in\textrm{ker}_{\varphi_{2}}$. Then by (\ref{equ:2.91})
and induction on $j_{q}$, we have $\mathcal{I} =\mathscr{A}$. Assume
that $\varphi_{q'} \equiv 0$, then $\mathcal{J}_{q'}=\mathbb{N}$.
Since $x^{-\sigma_{q}-\varepsilon,1_{[q']}}\in\mathcal{I}$ due to
$\varphi_{1}(-\sigma_{q}-\varepsilon)+2\not=0$, we have
\begin{equation}\label{equ:2.92}
\mathcal{I}\owns [x^{-\sigma_{q}-\varepsilon,1_{[q']}},
x^{\beta,\mathbf{j}}]=\epsilon_{q'}((j_{q'}
\varphi_{q}(\sigma_{q})+\varphi_{q}(\beta))x^{\beta,\mathbf{j}}+
j_{q}x^{\beta,\mathbf{j}-1_{[q]}})
\end{equation}
for $\beta\in\textrm{ker}_{\varphi_{2}}$. Then by (\ref{equ:2.92})
and induction on $j_{q}$, we have $\mathcal{I} =\mathscr{A}$.

Therefore, we have proved that $\mathcal{I} =\mathscr{A}$ if
$\boldsymbol{\mathcal{J}}\not=\{0\}$. So the first statement of
Theorem \ref{theorem:1} holds.\vspace{0.3cm}

\noindent{\bf \emph{Step 4. The second statement of Theorem
\ref{theorem:1} holds.}}\vspace{0.3cm}

If $\boldsymbol{\mathcal{J}}=\{0\}$, we have $\varphi_{p}\not\equiv
0$ for any $p\in\overline{1,2n}$. Fixing any $q\in\overline{3,2n}$,
we pick $ \tau \in \textrm{Rad}_{\phi} \cap(\cap _{q' \not= p \in
\overline{1,2n}} \textrm{ker}_{\varphi_p} )\backslash
\textrm{ker}_{\varphi_{q'}}$. Since $\varphi_{1}(\tau)+2\not=0$, we
have $x^{\tau, \mathbf{0}} \in \mathcal{I}$ by (\ref{equ:2.90}).
Then we get
\begin{equation}\label{equ:2.93}
\mathcal{I}\owns [x^{\tau,\mathbf{0}},
x^{\beta-\tau-\sigma_{q}-\varepsilon,\mathbf{0}}]=\epsilon_{q'}
\varphi_{q'}(\tau)\varphi_{q}(\beta-\sigma_{q})x^{\beta,\mathbf{0}}
\end{equation}
for $\beta\in\textrm{ker}_{\varphi_{2}}$. Thus we obtain
\begin{equation}\label{equ:2.94}
x^{\beta,\mathbf{0}}\in\mathcal{I}, \quad \textrm{for }
\beta\in\textrm{ker}_{\varphi_{2}} \textrm{ such that }
\varphi_{q}(\beta-\sigma_{q})\not=0.
\end{equation}
For $\beta \in \textrm{ker}_{\varphi_{2}}$ such that
$\varphi_{1}(\beta)+2=0$ and $\varphi_{q}(\beta-\sigma_{q})=0$ for
all $q\in\overline{3,2n}$, if $\beta\not\in\textrm{Rad}_{\phi}$,
there exists $\alpha\in\Gamma$ such that $\phi(\alpha,\beta)\not=0$.
Fix any such $\alpha$. If $\varphi_{1}(\alpha)+2=0$, replace
$\alpha$ by $\alpha+\varepsilon$. Then we can always assume that
$\varphi_{1}(\alpha)+2\not=0$. By (\ref{equ:2.88}) and
(\ref{equ:2.90}), we have $x^{\alpha,\mathbf{0}}\in\mathcal{I}$.
Then we get
\begin{eqnarray}
\mathcal{I} \owns
[x^{\alpha,\mathbf{0}},x^{\beta-\alpha-\varepsilon,\mathbf{0}}] & =&
\sum_{s=1}^{n-1}(\varphi_{2s+1}(\alpha)\varphi_{2s+2}(\beta)-\varphi_{2s+2}(\alpha)\varphi_{2s+1}(\beta))x^{\beta+\sigma_{2s+1},\mathbf{0}}\nonumber\\
& & +
\phi(\alpha,\beta)x^{\beta,\mathbf{0}}-(\varphi_{1}(\beta)+1)\varphi_{2}(\alpha)x^{\beta-\varepsilon,\mathbf{0}}.\label{equ:2.95}
\end{eqnarray}
Since $\varphi_{1}(\beta-\varepsilon)+2\not=0$, we have
$x^{\beta-\varepsilon,\mathbf{0}}\in\mathcal{I}$. Since
\begin{equation}
\varphi_{2s+1}(\beta+\sigma_{2s+1}-\sigma_{2s+1})=\varphi_{2s+1}(\sigma_{2s+1})\not=0
\end{equation}
or
\begin{equation}
\varphi_{2s+2}(\beta+\sigma_{2s+1}-\sigma_{2s+2})=\varphi_{2s+2}(\sigma_{2s+1})\not=0,
\end{equation}
we have $x^{\beta+\sigma_{2s+1},\mathbf{0}}\in\mathcal{I}$ for all
$s\in\overline{1,n-1}$. Then by (\ref{equ:2.95}) we get
\begin{equation}\label{equ:2.96}
x^{\beta,\mathbf{0}}\in\mathcal{I}\textrm{ for } \beta \in
\textrm{ker}_{\varphi_{2}}\backslash\textrm{Rad}_{\phi},
\varphi_{1}(\beta)+2=0,
\varphi_{q}(\beta-\sigma_{q})=0,q\in\overline{3,2n}.
\end{equation}
Thus we have obtained that
\begin{equation}\label{equ:2.97}
x^{\beta, \mathbf{0}}\in \mathcal{I}\quad \textrm{for}\ \beta\in
\Gamma \backslash \{\sigma\}
\end{equation}
by (\ref{equ:2.88}), (\ref{equ:2.90}), (\ref{equ:2.94}) and
(\ref{equ:2.96}). Set
\begin{equation}
\hat{\mathscr{H}}= \textrm{span} \{ x^{\alpha, \mathbf{0}} \mid
 \alpha \in \Gamma \backslash \{\sigma\}\}.
\end{equation}
Since $\boldsymbol{\mathcal{J}}=\{0\}$, by (\ref{equ:2.1}),
(\ref{equ:2.87}), (\ref{equ:2.89}), (\ref{equ:2.93}) and
(\ref{equ:2.95}), we have
\begin{equation}
\hat{\mathscr{H}}= [\mathscr{A}, \mathscr{A}].
\end{equation}
Replacing $\mathscr{A}$ by $\hat{\mathscr{H}}$ and $\Gamma$ by
$\Gamma\backslash \{\sigma\}$ in the above proof, we can also obtain
(\ref{equ:2.97}), which implies the simplicity of
\begin{displaymath}
\mathscr{H}^{(1)}=[\mathscr{H},
\mathscr{H}]=\hat{\mathscr{H}}/\mathbb{F}x^{-\varepsilon,
\mathbf{0}}.
\end{displaymath}
This completes the proof of Theorem \ref{theorem:1}.
$\qquad\qquad\qquad\bigbox$

\section{Proof of Theorem 2}
In this section, we will determine the irreducibility of the module
$\mathscr{A}_{\vec{\xi},f}$.

If there exists $\mu \in \Gamma$ such that (\ref{equ:1.10}) and
(\ref{equ:1.11}) hold, by (\ref{equ:1.9}) and (\ref{equ:2.1}) we
have the module isomorphism from $\mathscr{A}_{\vec{\xi},f}$ to
$\mathscr{A}_{\vec{\xi_{0}},0}$ in the following:
\begin{eqnarray}
\psi & : &  \mathscr{A}_{\vec{\xi},f}\longrightarrow  \mathscr{A}_{\vec{\xi_{0}},0},\\
 & & x^{\beta, \mathbf{j}}\longmapsto   x^{\beta+\mu-\varepsilon,
 \mathbf{j}}.
\end{eqnarray}
As a $\mathscr{H}$-module, $\mathscr{A}_{\vec{\xi_{0}},0}$ has a one
dimensional trivial submodule
$\mathbb{F}x^{-\varepsilon,\boldsymbol{0}}$. And
$\mathscr{A}_{\vec{\xi_{0}},0}/\mathbb{F}x^{-\varepsilon,\boldsymbol{0}}$
is irreducible when
$\boldsymbol{\mathcal{J}}\not=\{\boldsymbol{0}\}$ and indecomposable
when $\boldsymbol{\mathcal{J}}=\{\boldsymbol{0}\}$.

In the following, we assume that such $\mu$ does not exist.

Let $\mathcal {N}$ be any nonzero submodule of
$\mathscr{A}_{\vec{\xi},f}$. To prove the second statement of
Theorem \ref{theorem:2} is equivalent to proving that $\mathcal
{N}=\mathscr{A}_{\vec{\xi},f}$.

For $ k \in \mathbb{N} $, let
\begin{equation}
\mathscr{A}_{\vec{\xi},f}^{\langle k \rangle}=\textrm{span}
\{x^{\alpha,\mathbf{i}} \mid
(\alpha,\mathbf{i})\in\Gamma\times\boldsymbol{\mathcal{J}},|\mathbf{i}|\leq
k \}.
\end{equation}
Then we have
\begin{equation}
\mathscr{A}_{\vec{\xi},f}=\bigcup_{k=0}^{\infty}
\mathscr{A}_{\vec{\xi},f}^{\langle k \rangle}.
\end{equation}
Set
\begin{equation}\label{equ:3.1}
\hat{k'}=\textrm{min} \{ k \in \mathbb{N} \mid \mathcal{N}\cap
\mathscr{A}_{\vec{\xi},f}^{\langle k \rangle}\not = \{0\} \}.
\end{equation}
Then for any $u \in (\mathcal{N}\cap
\mathscr{A}_{\vec{\xi},f}^{\langle \hat{k'} \rangle}) \backslash
\{0\}$, we write it as
\begin{equation}\label{equ:3.2}
u= \sum_{\substack{(\alpha, \mathbf{i})\in
\Gamma \times \boldsymbol{\mathcal{J}},\\
|\mathbf{i}|= \hat{k'}}} b_{\alpha, \mathbf{i}} x^{\alpha,
\mathbf{i}}+u',
\end{equation}
where $b_{\alpha, \mathbf{i}} \in \mathbb{F}$ and $u' \in
\mathscr{A}_{\vec{\xi},f}^{\langle \hat{k'}-1 \rangle}$. Moreover,
we define
\begin{equation}\label{equ:3.3}
\kappa (u)=|\{\alpha\in\Gamma\mid b_{\alpha, \mathbf{i}} \not= 0 \
\textrm{for some}\ \mathbf{i}\in \boldsymbol{\mathcal{J}},\
|\mathbf{i}|= \hat{k'} \}|.
\end{equation}
By (\ref{equ:3.1}), $\kappa (u)>0$. Furthermore, we set
\begin{equation}\label{equ:3.4}
\kappa  = \textrm{min} \{\kappa  (w) \mid w \in (\mathcal{N}\cap
\mathscr{A}_{\vec{\xi},f}^{\langle \hat{k'} \rangle}) \backslash
\{0\} \}.
\end{equation}
Choose $u \in (\mathcal{N}\cap \mathscr{A}_{\vec{\xi},f}^{\langle
\hat{k'} \rangle}) \backslash \{0\} $ such that $\kappa (u) =
\kappa$. Write $u$ as in (\ref{equ:3.2}).

Quite similarly to the discussion in Step 1 and the second case of
Step 2 of the proof of Theorem 1, we can obtain some $x^{\alpha,
\mathbf{0}}\in \mathcal {N}$ with $\alpha\in\Gamma$. Here we omit
the details and let $\alpha$ be fixed. Then we will prove $\mathcal
{N}=\mathscr{A}_{\vec{\xi},f}$ by induction on $|\mathbf{i}|$.
Frequently, we will use the notation of Lie algebra homomorphism
\begin{equation}
\rho :  \mathscr{H} \longrightarrow
\textrm{GL}(\mathscr{A}_{\vec{\xi},f})
\end{equation}
to denote the module action in the rest of this paper.

Firstly we will prove $x^{\gamma, \mathbf{i}} \in \mathcal{N}$ for
all $\gamma\in\Gamma$ when $|\mathbf{i}|=0$. Namely, we will prove
$x^{\gamma, \mathbf{0}} \in \mathcal{N}$ for all $\gamma\in\Gamma$.
We give the proof in four cases.\vspace{0.2cm}

\noindent{\emph{Case 1. $\gamma\in\Gamma$ such that
\begin{equation}\label{equ:3.5}
(\varphi_{1}(\gamma)+1+\xi_{1})(\varphi_{2}(\alpha)+\xi_{2})-(\varphi_{2}(\gamma)+\xi_{2})(\varphi_{1}(\alpha)+\xi_{1})\not=0.
\end{equation}}}

By (\ref{equ:1.9}), we have
\begin{eqnarray}
 & & \mathcal{N}  \ni  \bar{x}^{\gamma-\alpha, \mathbf{0}}.x^{\alpha, \mathbf{0}}\nonumber\\
 & = &((\varphi_{1}(\gamma)+1+\xi_{1})(\varphi_{2}(\alpha)+\xi_{2})
-(\varphi_{2}(\gamma)+\xi_{2})(\varphi_{1}(\alpha)+\xi_{1}))x^{\gamma,\mathbf{0}}\nonumber\\
& &
+\sum_{s=1}^{n-1}\big((\varphi_{2s+1}(\gamma)+\xi_{2s+1})(\varphi_{2s+2}(\alpha)+\xi_{2s+2})
-(\varphi_{2s+2}(\gamma)+\xi_{2s+2})\nonumber\\
&
&\cdot(\varphi_{2s+1}(\alpha)+\xi_{2s+1})\big)x^{\gamma+\sigma_{2s+1}+\varepsilon,
\mathbf{0}}+(\phi(\gamma,\alpha)+f(\gamma-\alpha))x^{\gamma+\varepsilon,\mathbf{0}}.\label{equ:3.6}
\end{eqnarray}
If $\varphi_{2}\equiv 0$, we have $\mathcal{J}_{2}=\mathbb{N}$.
Moreover, if $\xi_{2}=0$, there does not exist any $\gamma\in\Gamma$
such that (\ref{equ:3.5}) holds. We assume that $\xi_{2}\not=0$.
Then we have
\begin{equation}\label{equ:3.7}
\bar{x}^{\varepsilon,\mathbf{0}}.(\bar{x}^{-\varepsilon,1_{[2]}}.x^{\beta,\mathbf{0}})=-2\xi_{2}(\varphi_{1}(\beta)+\xi_{1})x^{\beta,\mathbf{0}}.
\end{equation}
Repeatedly applying
$\rho(\bar{x}^{\varepsilon,\mathbf{0}})\rho(\bar{x}^{-\varepsilon,1_{[2]}})$
on $\mathcal{N}$, we obtain $x^{\gamma,\mathbf{0}} \in \mathcal{N}$
by Lemma \ref{lemma:1}, (\ref{equ:3.5}), (\ref{equ:3.6}) and
(\ref{equ:3.7}). If $\varphi_{2}\not\equiv 0$, picking $\tau \in
\textrm{Rad}_{\phi} \cap (\cap_{2\not=p \in \overline{1,2n}}
\textrm{ker}_{\varphi_{p}}) \backslash \textrm{ker}_{\varphi_{2}}$,
we have
\begin{equation}\label{equ:4.3}
\bar{x}^{\varepsilon-\tau,\mathbf{0}}.(\bar{x}^{\tau-\varepsilon,\mathbf{0}}.x^{\beta,\mathbf{0}})=S_{1}(\beta)x^{\beta,\mathbf{0}},
\end{equation}
where
$S_{1}(\beta)=-\varphi_{2}(\tau)(\varphi_{1}(\beta)+\xi_{1})(2(\varphi_{2}(\beta)+\xi_{2})+\varphi_{2}(\tau)(\varphi_{1}(\beta)+\xi_{1}+1))$.
Since we have
\begin{equation}
S_{1}(\gamma+\varepsilon)-S_{1}(\gamma)=-2\varphi_{2}(\tau)(\varphi_{2}(\tau)(\varphi_{1}(\gamma)+1+\xi_{1})+\varphi_{2}(\gamma)+\xi_{2}),
\end{equation}
and by (\ref{equ:3.5}) we have $\varphi_{1}(\gamma)+1+\xi_{1}\not=0$
or $\varphi_{2}(\gamma)+\xi_{2}\not=0$, we can restrict the
selection of $\tau$ so that
$S_{1}(\gamma+\varepsilon)-S_{1}(\gamma)\not=0$. Then repeatedly
applying
$\rho(\bar{x}^{\varepsilon-\tau,\mathbf{0}})\rho(\bar{x}^{\tau-\varepsilon,\mathbf{0}})$
on $\mathcal{N}$, we obtain that $x^{\gamma,\mathbf{0}} \in
\mathcal{N}$ by Lemma \ref{lemma:1}, (\ref{equ:3.5}),
(\ref{equ:3.6}) and (\ref{equ:4.3}).\vspace{0.2cm}

\noindent{\emph{Case 2. $\gamma\in\Gamma$ such that
\begin{eqnarray}
& &
(\varphi_{2t+1}(\gamma-\sigma_{2t+1})+\xi_{2t+1})(\varphi_{2t+2}(\alpha)+\xi_{2t+2})\nonumber\\
&  & -
(\varphi_{2t+2}(\gamma-\sigma_{2t+1})+\xi_{2t+2})(\varphi_{2t+1}(\alpha)+\xi_{2t+1})\not=0\label{equ:3.9}
\end{eqnarray}
for some $t\in\overline{1,n-1}$.}}

By (\ref{equ:1.9}), we have
\begin{eqnarray}
& &  \mathcal{N}  \ni  \bar{x}^{\gamma-\sigma_{2t+1}-\varepsilon-\alpha, \mathbf{0}}.x^{\alpha, \mathbf{0}}\nonumber\\
 & = & ((\varphi_{1}(\gamma)+\xi_{1})(\varphi_{2}(\alpha)+\xi_{2})
-(\varphi_{2}(\gamma)+\xi_{2})(\varphi_{1}(\alpha)+\xi_{1}))x^{\gamma-\sigma_{2t+1}-\varepsilon,\mathbf{0}}\nonumber\\
& & +\sum_{\substack{s=1 \\ s\not=t}}^{n-1} \big(
(\varphi_{2s+1}(\gamma)+\xi_{2s+1})(\varphi_{2s+2}(\alpha)+\xi_{2s+2})
-(\varphi_{2s+2}(\gamma)+\xi_{2s+2})\nonumber\\
& &\cdot(\varphi_{2s+1}(\alpha)+\xi_{2s+1})
\big)x^{\gamma+\sigma_{2s+1}-\sigma_{2t+1}, \mathbf{0}}+(\phi(\gamma,\alpha)+f(\gamma-\alpha))x^{\gamma-\sigma_{2t+1},\mathbf{0}}\nonumber\\
& &+ \big(
(\varphi_{2t+1}(\gamma-\sigma_{2t+1})+\xi_{2t+1})(\varphi_{2t+2}(\alpha)+\xi_{2t+2})-
(\varphi_{2t+1}(\alpha)+\xi_{2t+1})\nonumber\\
&  & \cdot (\varphi_{2t+2}(\gamma-\sigma_{2t+1})+\xi_{2t+2})
\big)x^{\gamma, \mathbf{0}} . \label{equ:3.10}
\end{eqnarray}
If $\sigma_{2t+1}=0$, we have $\varphi_{2t+1}\equiv
\varphi_{2t+2}\equiv 0$ by (\ref{equ:1.4}). Then there does not
exist any $\gamma\in\Gamma$ such that (\ref{equ:3.9}) holds. Assume
that $\sigma_{2t+1}\not=0$ , we have
$\varphi_{2t+1}(\sigma_{2t+1})\not=0$ or
$\varphi_{2t+2}(\sigma_{2t+1})\not=0$ by (\ref{equ:1.4}). Choose
$q\in\{2t+1,2t+2\}$ such that $\varphi_{q}(\sigma_{2t+1})\not=0$.
Then $\varphi_{q}\not\equiv 0$, and we have that
$\sigma_{q}=\sigma_{q'}=\sigma_{2t+1}$.

If $\varphi_{q'}\equiv 0$, we have $\mathcal{J}_{q'}=\mathbb{N}$.
Moreover, if $\xi_{q'}=0$, there does not exist any
$\gamma\in\Gamma$ such that (\ref{equ:3.9}) holds. Suppose that
$\xi_{q'}\not=0$. Then we have
\begin{equation}\label{equ:3.11}
\bar{x}^{-2\sigma_{q}-\varepsilon,\mathbf{0}}.(\bar{x}^{-\varepsilon,1_{[q']}}.x^{\beta,\mathbf{0}})=2\xi_{q'}\varphi_{q}(\sigma_{q})(\varphi_{q}(\beta)+\xi_{q})x^{\beta,\mathbf{0}}.
\end{equation}
Repeatedly applying
$\rho(\bar{x}^{-2\sigma_{q}-\varepsilon,\mathbf{0}})\rho(\bar{x}^{-\varepsilon,1_{[q']}})$
on $\mathcal{N}$, we obtain that $x^{\gamma,\mathbf{0}} \in
\mathcal{N}$ by Lemma \ref{lemma:1}, (\ref{equ:3.9}),
(\ref{equ:3.10}) and (\ref{equ:3.11}).

If $\varphi_{q'}\not\equiv 0$, picking $\tau' \in
\textrm{Rad}_{\phi} \cap (\cap_{q'\not=p \in \overline{1,2n}}
\textrm{ker}_{\varphi_{p}}) \backslash \textrm{ker}_{\varphi_{q'}}$,
we have
\begin{equation}\label{equ:4.4}
\bar{x}^{-\tau'-2\sigma_{q}-\varepsilon,\mathbf{0}}.(\bar{x}^{\tau'-\varepsilon,\mathbf{0}}.x^{\beta,\mathbf{0}})=S_{2}(\beta)x^{\beta,\mathbf{0}},
\end{equation}
where
\begin{eqnarray}
S_{2}(\beta) &=&
\varphi_{q'}(\tau')(\varphi_{q}(\beta)+\xi_{q})\big(2\varphi_{q}(\sigma_{q})(\varphi_{q'}(\beta+\tau'+\sigma_{q})+\xi_{q'})\nonumber\\
& &
-\varphi_{q'}(\tau'+2\sigma_{q})(\varphi_{q}(\beta+\sigma_{q})+\xi_{q})\big).
\end{eqnarray}
Since
\begin{eqnarray}
S_{2}(\gamma)-S_{2}(\gamma-\sigma_{q})  &=&
2\varphi_{q'}(\tau')\varphi_{q}(\sigma_{q})\big(\varphi_{q}(\sigma_{q})(\varphi_{q'}(\gamma-\sigma_{q})+\xi_{q'})\nonumber\\
& &
-\varphi_{q'}(\tau'+\sigma_{q})(\varphi_{q}(\gamma-\sigma_{q})+\xi_{q})\big),
\end{eqnarray}
and by (\ref{equ:3.9}) we have
$\varphi_{q}(\gamma-\sigma_{q})+\xi_{q}\not=0$ or
$\varphi_{q'}(\gamma-\sigma_{q})+\xi_{q'}\not=0$, we can restrict
the selection of $\tau'$ so that
$S_{2}(\gamma)-S_{2}(\gamma-\sigma_{q})\not=0$. Then repeatedly
applying
$\rho(\bar{x}^{-\tau'-2\sigma_{q}-\varepsilon,\mathbf{0}})\rho(\bar{x}^{\tau'-\varepsilon,\mathbf{0}})$
on $\mathcal{N}$, we obtain that $x^{\gamma,\mathbf{0}} \in
\mathcal{N}$ by Lemma \ref{lemma:1}, (\ref{equ:3.9}),
(\ref{equ:3.10}) and (\ref{equ:4.4}).\vspace{0.2cm}

\noindent{\emph{Case 3. $\gamma\in\Gamma$ such that
\begin{equation}\label{equ:3.13}
\phi(\gamma, \alpha)+f(\gamma-\alpha)\not=0.
\end{equation}}}

By (\ref{equ:1.9}), we have
\begin{eqnarray}
 & & \mathcal{N}  \ni  \bar{x}^{\gamma-\varepsilon-\alpha, \mathbf{0}}.x^{\alpha, \mathbf{0}}\nonumber\\
 & = &\big((\varphi_{1}(\gamma)+\xi_{1})(\varphi_{2}(\alpha)+\xi_{2})
-(\varphi_{2}(\gamma)+\xi_{2})(\varphi_{1}(\alpha)+\xi_{1})\big)x^{\gamma-\varepsilon,\mathbf{0}}\nonumber\\
& & + \sum_{s=1}^{n-1}
\big((\varphi_{2s+1}(\gamma)+\xi_{2s+1})(\varphi_{2s+2}(\alpha)+\xi_{2s+2})
-(\varphi_{2s+2}(\gamma)+\xi_{2s+2})\nonumber\\
& &\cdot
(\varphi_{2s+1}(\alpha)+\xi_{2s+1})\big)x^{\gamma+\sigma_{2s+1},
\mathbf{0}}
+(\phi(\gamma,\alpha)+f(\gamma-\alpha))x^{\gamma,\mathbf{0}}.\label{equ:3.14}
\end{eqnarray}
Similarly to the discussion from (\ref{equ:3.6}) to (\ref{equ:4.4}),
we can obtain that $x^{\gamma,\mathbf{0}} \in
\mathcal{N}$.\vspace{0.2cm}

\noindent{\emph{Case 4. $\gamma\in\Gamma$ such that
\begin{equation}\label{equ:3.15}
(\varphi_{1}(\gamma)+1+\xi_{1})(\varphi_{2}(\alpha)+\xi_{2})-(\varphi_{2}(\gamma)+\xi_{2})(\varphi_{1}(\alpha)+\xi_{1})=0,
\end{equation}
\begin{eqnarray}
& &
(\varphi_{2t+1}(\gamma-\sigma_{2t+1})+\xi_{2t+1})(\varphi_{2t+2}(\alpha)+\xi_{2t+2})\nonumber\\
&  & -
(\varphi_{2t+2}(\gamma-\sigma_{2t+1})+\xi_{2t+2})(\varphi_{2t+1}(\alpha)+\xi_{2t+1})=0\label{equ:3.16}
\end{eqnarray}
for all $t\in\overline{1,n-1}$, and
\begin{equation}\label{equ:3.17}
\phi(\gamma, \alpha)+f(\gamma-\alpha)=0.
\end{equation}}}

\noindent{\emph{Subcase 4.1.
$\varphi_{2}(\gamma)+\xi_{2}\not=0$.}}\vspace{0.1cm}

Choose $\tau \in \textrm{Rad}_{\phi} \cap (\cap_{p \in
\overline{2,2n}} \textrm{ker}_{\varphi_{p}}) \backslash
\textrm{ker}_{\varphi_{1}}$ such that $\varphi_{1}(\tau)+1\not=0$.

If $\varphi_{2}(\alpha)+\xi_{2}\not=0$, we have
\begin{equation}\label{equ:3.18}
\mathcal{N} \ni
\bar{x}^{\tau,\mathbf{0}}.x^{\alpha,\mathbf{0}}=(\varphi_{1}(\tau)+1)(\varphi_{2}(\alpha)+\xi_{2})x^{\alpha+\tau,\mathbf{0}}.
\end{equation}
Then we go back to \emph{Case 1} with $x^{\alpha,\mathbf{0}}$
replaced by $x^{\alpha+\tau,\mathbf{0}}$.

If $\varphi_{2}(\alpha)+\xi_{2}=0$, by (\ref{equ:3.15}) we have
$\varphi_{1}(\alpha)+\xi_{1}=0$.

Moreover, if there exists $q\in\overline{3,2n}$ such that
$\varphi_{q}(\alpha)+\xi_{q}\not=0$ and $\varphi_{q'}\not\equiv0$,
picking $\tau' \in \textrm{Rad}_{\phi} \cap (\cap_{q'\not=p \in
\overline{1,2n}} \textrm{ker}_{\varphi_{p}}) \backslash
\textrm{ker}_{\varphi_{q'}}$, we have
\begin{equation}\label{equ:3.19}
\mathcal{N} \ni
\bar{x}^{\tau+\tau',\mathbf{0}}.x^{\alpha,\mathbf{0}}=\epsilon_{q'}\varphi_{q'}(\tau')(\varphi_{q}(\alpha)+\xi_{q})x^{\alpha+\tau+\tau'+\sigma_{q}+\varepsilon,\mathbf{0}}.
\end{equation}
Then we go back to \emph{Case 1} with $x^{\alpha,\mathbf{0}}$
replaced by
$x^{\alpha+\tau+\tau'+\sigma_{q}+\varepsilon,\mathbf{0}}$.

If there exists $q\in\overline{3,2n}$ such that
$\varphi_{q}(\alpha)+\xi_{q}\not=0$ and $\varphi_{q'}\equiv0$, we
have
\begin{equation}\label{equ:3.20}
\mathcal{N} \ni
\bar{x}^{\tau,1_{[q']}}.x^{\alpha,\mathbf{0}}=\epsilon_{q'}(\varphi_{q}(\alpha)+\xi_{q})x^{\alpha+\tau+\sigma_{q}+\varepsilon,\mathbf{0}}.
\end{equation}
We go back to \emph{Case 1} with $x^{\alpha,\mathbf{0}}$ replaced by
$x^{\alpha+\tau+\sigma_{q}+\varepsilon,\mathbf{0}}$.

If $\varphi_{q}(\alpha)+\xi_{q}=0$ for all $q\in\overline{3,2n}$,
then there exists $\theta\in\Gamma$ such that
\begin{equation}
\phi(\theta, \alpha)+f(\theta)\not=0.
\end{equation}
We rechoose $\tau \in \textrm{Rad}_{\phi} \cap (\cap_{p \in
\overline{2,2n}} \textrm{ker}_{\varphi_{p}}) \backslash
\textrm{ker}_{\varphi_{1}}$ such that
$\varphi_{1}(\alpha+\theta+\tau+\varepsilon)+\xi_{1}\not=0$. Then we
have
\begin{equation}\label{equ:3.21}
\mathcal{N} \ni
\bar{x}^{\theta+\tau,\mathbf{0}}.x^{\alpha,\mathbf{0}}=(\phi(\theta,\alpha)+f(\theta))x^{\alpha+\theta+\tau+\varepsilon,\mathbf{0}}.
\end{equation}
We go back to \emph{Case 1} or the first condition of \emph{Subcase
4.1} with $x^{\alpha,\mathbf{0}}$ replaced by
$x^{\alpha+\theta+\tau+\varepsilon,\mathbf{0}}$. So we obtain that
$x^{\gamma,\mathbf{0}} \in \mathcal{N}$.\vspace{0.1cm}

\noindent{\emph{Subcase 4.2. $\varphi_{2}(\gamma)+\xi_{2}=0$ and
$\varphi_{1}(\gamma)+1+\xi_{1}\not=0$.}}\vspace{0.1cm}

By (\ref{equ:3.15}) we have that $\varphi_{2}(\alpha)+\xi_{2}=0$. If
$\varphi_{2}\not\equiv 0$, we can prove $x^{\gamma,\mathbf{0}} \in
\mathcal{N}$ similarly to \emph{Subcase 4.1}. Here we omit the
details. We assume that $\varphi_{2}\equiv 0$. Then we get
$\mathcal{J}_{2}=\mathbb{N}$ and $\xi_{2}=0$.

If $\varphi_{1}(\alpha)+\xi_{1}\not=0$, we have
\begin{eqnarray}
& & \mathcal{N}  \ni
\bar{x}^{\gamma-\alpha,1_{[2]}}.x^{\alpha,\mathbf{0}}\nonumber\\
& = & \sum_{s=1}^{n-1} \big(
(\varphi_{2s+1}(\gamma)+\xi_{2s+1})(\varphi_{2s+2}(\alpha)+\xi_{2s+2})
-(\varphi_{2s+2}(\gamma)+\xi_{2s+2})\nonumber\\
& &
\cdot(\varphi_{2s+1}(\alpha)+\xi_{2s+1})\big)x^{\gamma+\sigma_{2s+1}+\varepsilon,1_{[2]}}
-(\varphi_{1}(\alpha)+\xi_{1})x^{\gamma,\mathbf{0}}\label{equ:3.22}.
\end{eqnarray}
Since we have
\begin{eqnarray}
& &
\bar{x}^{\varepsilon,1_{[2]}}.(\bar{x}^{-\varepsilon,1_{[2]}}.x^{\gamma+\sigma_{2s+1}+\varepsilon,1_{[2]}})\nonumber\\
& = &
(\varphi_{1}(\gamma)+\xi_{1}+1)(\varphi_{1}(\gamma)+\xi_{1}-2)x^{\gamma+\sigma_{2s+1}+\varepsilon,1_{[2]}}\label{equ:3.23}
\end{eqnarray}
for $s\in\overline{1,n-1}$, and
\begin{equation}\label{equ:3.24}
\bar{x}^{\varepsilon,1_{[2]}}.(\bar{x}^{-\varepsilon,1_{[2]}}.x^{\gamma,\mathbf{0}})=(\varphi_{1}(\gamma)+\xi_{1})(\varphi_{1}(\gamma)+\xi_{1}-1)x^{\gamma,\mathbf{0}}.
\end{equation}
Repeatedly applying
$\rho(\bar{x}^{\varepsilon,1_{[2]}})\rho(\bar{x}^{-\varepsilon,1_{[2]}})$
on $\mathcal{N}$, we obtain $x^{\gamma,\mathbf{0}} \in \mathcal{N}$
by Lemma \ref{lemma:1}, (\ref{equ:3.22}), (\ref{equ:3.23}) and
(\ref{equ:3.24}).

If $\varphi_{1}(\alpha)+\xi_{1}=0$, pick $\tau \in
\textrm{Rad}_{\phi} \cap (\cap_{p \in \overline{2,2n}}
\textrm{ker}_{\varphi_{p}}) \backslash \textrm{ker}_{\varphi_{1}}$
such that $\varphi_{1}(\tau)+1\not=0$.

Moreover, if there exists $q\in\overline{3,2n}$ such that
$\varphi_{q}(\alpha)+\xi_{q}\not=0$ and $\varphi_{q'}\not\equiv0$,
picking $\tau' \in \textrm{Rad}_{\phi} \cap (\cap_{q'\not=p \in
\overline{1,2n}} \textrm{ker}_{\varphi_{p}}) \backslash
\textrm{ker}_{\varphi_{q'}}$, we have
\begin{equation}\label{equ:3.25}
\mathcal{N} \ni
\bar{x}^{\tau+\tau',\mathbf{0}}.x^{\alpha,\mathbf{0}}=\epsilon_{q'}\varphi_{q'}(\tau')(\varphi_{q}(\alpha)+\xi_{q})x^{\alpha+\tau+\tau'+\sigma_{q}+\varepsilon,\mathbf{0}}.
\end{equation}
Since
$\varphi_{1}(\alpha+\tau+\tau'+\sigma_{q}+\varepsilon)+\xi_{1}\not=0$,
we go back to the previous condition from (\ref{equ:3.22}) to
(\ref{equ:3.24}) with $x^{\alpha,\mathbf{0}}$ replaced by
$x^{\alpha+\tau+\tau'+\sigma_{q}+\varepsilon,\mathbf{0}}$. So we
obtain that $x^{\gamma,\mathbf{0}} \in \mathcal{N}$.

If there exists $q\in\overline{3,2n}$ such that
$\varphi_{q}(\alpha)+\xi_{q}\not=0$ and $\varphi_{q'}\equiv0$, then
we have
\begin{equation}\label{equ:3.26}
\mathcal{N} \ni
\bar{x}^{\tau,1_{[q']}}.x^{\alpha,\mathbf{0}}=\epsilon_{q'}(\varphi_{q}(\alpha)+\xi_{q})x^{\alpha+\tau+\sigma_{q}+\varepsilon,\mathbf{0}}.
\end{equation}
Since
$\varphi_{1}(\alpha+\tau+\sigma_{q}+\varepsilon)+\xi_{1}\not=0$, we
go back to the condition from (\ref{equ:3.22}) to (\ref{equ:3.24})
with $x^{\alpha,\mathbf{0}}$ replaced by
$x^{\alpha+\tau+\sigma_{q}+\varepsilon,\mathbf{0}}$. So we obtain
that $x^{\gamma,\mathbf{0}} \in \mathcal{N}$.

If $\varphi_{q}(\alpha)+\xi_{q}=0$ for all $q\in\overline{3,2n}$,
then there exists $\theta\in\Gamma$ such that
\begin{equation}
\phi(\theta, \alpha)+f(\theta)\not=0.
\end{equation}
We rechoose $\tau \in \textrm{Rad}_{\phi} \cap (\cap_{p \in
\overline{2,2n}} \textrm{ker}_{\varphi_{p}}) \backslash
\textrm{ker}_{\varphi_{1}}$ such that
$\varphi_{1}(\alpha+\theta+\tau+\varepsilon)+\xi_{1}\not=0$, then we
have
\begin{equation}\label{equ:3.27}
\mathcal{N} \ni
\bar{x}^{\theta+\tau,\mathbf{0}}.x^{\alpha,\mathbf{0}}=(\phi(\theta,\alpha)+f(\theta))x^{\alpha+\theta+\tau+\varepsilon,\mathbf{0}}.
\end{equation}
We go back to the condition from (\ref{equ:3.22}) to
(\ref{equ:3.24}) with $x^{\alpha,\mathbf{0}}$ replaced by
$x^{\alpha+\theta+\tau+\varepsilon,\mathbf{0}}$. So we obtain that
$x^{\gamma,\mathbf{0}} \in \mathcal{N}$.\vspace{0.1cm}

\noindent{\emph{Subcase 4.3. $\varphi_{2}(\gamma)+\xi_{2}=0$,
$\varphi_{1}(\gamma)+1+\xi_{1}=0$ and there exists
$q\in\overline{3,2n}$ such that
$\varphi_{q}(\gamma-\sigma_{q})+\xi_{q}\not=0$.}}\vspace{0.1cm}

The proof for $x^{\gamma,\mathbf{0}} \in \mathcal{N}$ in this
subcase is similar to the proof in \emph{Subcase 4.2}. Here we omit
the details.\vspace{0.1cm}

\noindent{\emph{Subcase 4.4. $\varphi_{2}(\gamma)+\xi_{2}=0$,
$\varphi_{1}(\gamma)+1+\xi_{1}=0$ and
$\varphi_{q}(\gamma-\sigma_{q})+\xi_{q}=0$ for all
$q\in\overline{3,2n}$.}}\vspace{0.1cm}

Since we assumed that $\mu$ defined in Theorem \ref{theorem:2} does
not exist, there exists $\theta \in \Gamma$ such that
\begin{equation}
\phi(\theta, \gamma+\varepsilon-\sum_{s=1}^{n-1}\sigma_{2s+1}) +
f(\theta)\not=0.
\end{equation}
Namely, we have $\phi(\gamma,\theta) - f(\theta)\not=0$.

If $\varphi_{2}(\alpha)+\xi_{2}\not=0$, choose $\tau \in
\textrm{Rad}_{\phi} \cap (\cap_{p \in \overline{2,2n}}
\textrm{ker}_{\varphi_{p}}) \backslash \textrm{ker}_{\varphi_{1}}$
such that
\begin{equation}
(\varphi_{1}(\tau+\theta)+1)(\varphi_{2}(\alpha)+\xi_{2})-\varphi_{2}(\theta)(\varphi_{1}(\alpha)+\xi_{1})\not=0
\end{equation}
and
\begin{equation}
\varphi_{1}(\alpha+\tau+\theta)+1+\xi_{1}\not=0.
\end{equation}
Then we have
\begin{eqnarray}
 & & \mathcal{N}  \ni  \bar{x}^{\tau+\theta, \mathbf{0}}.x^{\alpha, \mathbf{0}}\nonumber\\
 & = & ((\varphi_{1}(\tau+\theta)+1)(\varphi_{2}(\alpha)+\xi_{2})
-\varphi_{2}(\theta)(\varphi_{1}(\alpha)+\xi_{1}))x^{\alpha+\tau+\theta,\mathbf{0}}
\nonumber\\
& &
+(\phi(\theta,\alpha)+f(\theta))x^{\alpha+\tau+\theta+\varepsilon,\mathbf{0}}
+\sum_{s=1}^{n-1}
\big(\varphi_{2s+1}(\theta)(\varphi_{2s+2}(\alpha)+\xi_{2s+2})
\nonumber\\
& &
-\varphi_{2s+2}(\theta)(\varphi_{2s+1}(\alpha)+\xi_{2s+1})\big)x^{\alpha+\tau+\theta+\sigma_{2s+1}+\varepsilon,
\mathbf{0}}.\label{equ:3.28}
\end{eqnarray}
Similarly to the discussion from (\ref{equ:3.6}) to (\ref{equ:4.3})
in \emph{Case 1}, we obtain that $x^{\alpha+\tau+\theta,\mathbf{0}}
\in \mathcal{N}$. Then we go back to \emph{Case 3} with
$x^{\alpha,\mathbf{0}}$ replaced by
$x^{\alpha+\tau+\theta,\mathbf{0}}$.

Suppose $\varphi_{2}(\alpha)+\xi_{2}=0$ and
$\varphi_{1}(\alpha)+\xi_{1}\not=0$. If $\varphi_{2}\not\equiv 0$,
choose
\begin{equation}
\tau \in \textrm{Rad}_{\phi} \cap (\cap_{2\not=p \in
\overline{1,2n}} \textrm{ker}_{\varphi_{p}}) \backslash
\textrm{ker}_{\varphi_{2}}
\end{equation}
such that
\begin{equation}
\varphi_{2}(\tau+\theta)\not=0 \textrm{ and }
\varphi_{2}(\alpha+\tau+\theta)+\xi_{2}\not=0.
\end{equation}
Then we have
\begin{eqnarray}
 & & \mathcal{N}  \ni  \bar{x}^{\tau+\theta, \mathbf{0}}.x^{\alpha, \mathbf{0}}\nonumber\\
 & = & (\phi(\theta,\alpha)+f(\theta))x^{\alpha+\tau+\theta+\varepsilon,\mathbf{0}}
-\varphi_{2}(\tau+\theta)(\varphi_{1}(\alpha)+\xi_{1})x^{\alpha+\tau+\theta,\mathbf{0}}
\nonumber\\
& &  +\sum_{s=1}^{n-1}
\big(\varphi_{2s+1}(\theta)(\varphi_{2s+2}(\alpha)+\xi_{2s+2})
-\varphi_{2s+2}(\theta)(\varphi_{2s+1}(\alpha)+\xi_{2s+1})\big)\nonumber\\
& &\cdot x^{\alpha+\tau+\theta+\sigma_{2s+1}+\varepsilon,
\mathbf{0}}.\label{equ:3.29}
\end{eqnarray}
Similarly to the discussion from (\ref{equ:3.6}) to (\ref{equ:4.3})
in \emph{Case 1}, we obtain that $x^{\alpha+\tau+\theta,\mathbf{0}}
\in \mathcal{N}$. Then we go back to \emph{Case 3} with
$x^{\alpha,\mathbf{0}}$ replaced by
$x^{\alpha+\tau+\theta,\mathbf{0}}$. Assume that $\varphi_{2}\equiv
0$, then we have $\mathcal{J}_{2}=\mathbb{N}$ and $\xi_{2}=0$. Since
we have
\begin{eqnarray}
 & & \mathcal{N}  \ni  \bar{x}^{\theta,1_{[2]}}.x^{\alpha, \mathbf{0}}\nonumber\\
 & = &
\sum_{s=1}^{n-1}
(\varphi_{2s+1}(\theta)(\varphi_{2s+2}(\alpha)+\xi_{2s+2})
-\varphi_{2s+2}(\theta)(\varphi_{2s+1}(\alpha)+\xi_{2s+1}))x^{\alpha+\theta+\sigma_{2s+1}+\varepsilon,
1_{[2]}}\nonumber\\
& &
+(\phi(\theta,\alpha)+f(\theta))x^{\alpha+\theta+\varepsilon,1_{[2]}}
-(\varphi_{1}(\alpha)+\xi_{1})x^{\alpha+\theta,\mathbf{0}}.\label{equ:3.30}
\end{eqnarray}
Similarly to the discussion from (\ref{equ:3.22}) to
(\ref{equ:3.24}), we obtain that $x^{\alpha+\theta,\mathbf{0}} \in
\mathcal{N}$. We go back to \emph{Case 3} with
$x^{\alpha,\mathbf{0}}$ replaced by $x^{\alpha+\theta,\mathbf{0}}$,
then we obtain $x^{\gamma,\mathbf{0}} \in \mathcal{N}$.

Next we consider the case that $\varphi_{2}(\alpha)+\xi_{2}=0$,
$\varphi_{1}(\alpha)+\xi_{1}=0$ and there exists
$q\in\overline{3,2n}$ such that $\varphi_{q}(\alpha)+\xi_{q}\not=0$.
Choose $\tau \in \textrm{Rad}_{\phi} \cap (\cap_{p \in
\overline{2,2n}} \textrm{ker}_{\varphi_{p}}) \backslash
\textrm{ker}_{\varphi_{1}}$ such that $\varphi_{1}(\tau)+1\not=0$.
If $\varphi_{q'}\not\equiv0$, picking $\tau' \in \textrm{Rad}_{\phi}
\cap (\cap_{q'\not=p \in \overline{1,2n}}
\textrm{ker}_{\varphi_{p}}) \backslash \textrm{ker}_{\varphi_{q'}}$,
we have
\begin{equation}\label{equ:4.5}
\mathcal{N}\ni \bar{x}^{\tau+\tau',\mathbf{0}}.x^{\alpha,
\mathbf{0}}=\epsilon_{q'}\varphi_{q'}(\tau')(\varphi_{q}(\alpha)+\xi_{q})x^{\alpha+\tau+\tau'+\sigma_{q}+\varepsilon,
\mathbf{0}}.
\end{equation}
Since
$\varphi_{2}(\alpha+\tau+\tau'+\sigma_{q}+\varepsilon)+\xi_{2}=0$
and
$\varphi_{1}(\alpha+\tau+\tau'+\sigma_{q}+\varepsilon)+\xi_{1}\not=0$,
we go back to the condition from (\ref{equ:3.29}) to
(\ref{equ:3.30}) with $x^{\alpha,\mathbf{0}}$ replaced by
$x^{\alpha+\tau+\tau'+\sigma_{q}+\varepsilon, \mathbf{0}}$. If
$\varphi_{q'}\equiv0$, we have $\mathcal{J}_{q'}=\mathbb{N}$. Then
we get
\begin{equation}\label{equ:4.6}
\mathcal{N}\ni \bar{x}^{\tau,1_{[q']}}.x^{\alpha,
\mathbf{0}}=\epsilon_{q'}(\varphi_{q}(\alpha)+\xi_{q})x^{\alpha+\tau+\sigma_{q}+\varepsilon,
\mathbf{0}}.
\end{equation}
Since $\varphi_{2}(\alpha+\tau+\sigma_{q}+\varepsilon)+\xi_{2}=0$
and $\varphi_{1}(\alpha+\tau+\sigma_{q}+\varepsilon)+\xi_{1}\not=0$,
we go back to the condition from (\ref{equ:3.29}) to
(\ref{equ:3.30}) with $x^{\alpha,\mathbf{0}}$ replaced by
$x^{\alpha+\tau+\sigma_{q}+\varepsilon, \mathbf{0}}$.

If $\varphi_{q}(\alpha)+\xi_{q}=0$ for all $q\in\overline{1,2n}$,
then there exists $\theta' \in \Gamma$ such that
\begin{equation}
\phi(\theta',\alpha)+f(\theta')\not=0.
\end{equation}
Choose $m\in\mathbb{N}$ such that
\begin{equation}
\phi(m\theta'+\theta, \alpha)+f(m\theta'+\theta)\not =0
\end{equation}
and
\begin{equation}
\phi(\gamma,m\theta'+\theta)-f(m\theta'+\theta)\not =0.
\end{equation}
Then we have
\begin{equation}\label{equ:3.31}
\mathcal{N} \ni \bar{x}^{m\theta'+\theta,\mathbf{0}}.x^{\alpha,
\mathbf{0}}=(\phi(m\theta'+\theta,
\alpha)+f(m\theta'+\theta))x^{\alpha+m\theta'+\theta+\varepsilon,\mathbf{0}}.
\end{equation}
We go back to \emph{Case 3} with $x^{\alpha,\mathbf{0}}$ replaced by
$x^{\alpha+m\theta'+\theta+\varepsilon,\mathbf{0}}$, then we obtain
that $x^{\gamma,\mathbf{0}} \in \mathcal{N}$.

Therefore we have proved that $x^{\gamma, \mathbf{i}} \in
\mathcal{N}$ for all $\gamma\in\Gamma$ when $|\mathbf{i}|=0$. Fixing
$0<k\in\mathbb{N}$, assume that $x^{\gamma, \mathbf{i}} \in
\mathcal{N}$ for all $\gamma\in\Gamma$ with $|\mathbf{i}|<k$. Then
we need to prove that $x^{\gamma, \mathbf{i}} \in \mathcal{N}$ for
all $\gamma\in\Gamma$ with $|\mathbf{i}|=k$. The proof is similar to
the proof for $|\mathbf{i}|=0$. Here we omit the details.

This completes the proof of Theorem \ref{theorem:2}.
$\qquad\qquad\qquad\bigbox$

\section*{Acknowledgment}
I would like to express my deep gratitude to Professor Xiaoping Xu
for all his advice, instructions and encouragements.

\end{document}